\documentclass[a4paper,12pt]{article} 
\newdimen\paperhight
\usepackage{amsmath}
\usepackage{amssymb}
\usepackage{amsfonts}
\usepackage[usenames]{color}

\newcommand{\dsp}{\displaystyle}
\newcommand{\dto}{\downarrow}

\newcommand{\pd}{\partial} 

\newcommand{\pr}{\par \vspace{3mm}\noindent [{\bf Proof}] \qquad}
\newcommand{\prend}{\hfill \qed \par \vspace{3mm}}

\newcommand{\qed}{\quad\hbox{\rule[-2pt]{3pt}{6pt}}\par\vspace{3mm}}

\newcommand{\1}{{\bf 1}} 
 
\newcommand{\C}{\mathbb C} 
\newcommand{\Z}{\mathbb Z}

\newcommand{\N}{\mathbb N}

\newcommand{\CI}{{\cal I}}
\newcommand{\CF}{{\cal F}}
\newcommand{\CL}{{\cal L}}

\newcommand{\CU}{{\cal U}}

\newcommand{\CY}{{\cal Y}}

\newcommand{\CH}{{\cal H}}

\newcommand{\wt}{{\rm wt}}
\newcommand{\Hom}{{\rm Hom}}

\newcommand{\End}{{\rm End}}
\newcommand{\Image}{{\rm Im}}
\newcommand{\Ker}{{\rm Ker}}
\newcommand{\Res}{{\rm Res}}
\newcommand{\Tr}{{\rm Tr}}
\newtheorem{thm}{Theorem}
\newtheorem{prn}[thm]{Proposition}
\newtheorem{dfn}[thm]{Definition}
\newtheorem{cry}[thm]{Corollary}
\newtheorem{lmm}[thm]{Lemma}

\begin{document}
\title{Flatness and Semi-Rigidity of Vertex Operator Algebras}
\author{\begin{tabular}{c}
Masahiko Miyamoto\\
Institute of Mathematics, \\
University of Tsukuba, \\
Tsukuba, 305 Japan \end{tabular}}
\date{}
\maketitle

\begin{abstract}
In VOA theories, most of the general theorems are proved under the assumption of 
rationality and $C_2$-cofiniteness. 
In this paper, we obtain several general theorems without the assumption of 
rationality so that we can use them for proving rationality of given $C_2$-cofinite VOAs. 
For example, we apply them to orbifold models and show that 
if $g$ is a finite automorphism of a rational $C_2$-cofinite VOA $T$ of 
CFT-type with $T'\cong T$ and 
a fixed point subVOA $V:=T^g$ is a $C_2$-cofinite, then $V$ is also rational. 
We also investigate several results under weaker conditions. 
\end{abstract} 

\section{Introduction}
In this paper, we summarize 
results given in three preprints \cite{M2},\cite{M3} and \cite{M4} and extend some of them 
under weaker conditions.   
In the research of vertex operator algebras (shortly VOAs), 
most of the general theorems are proved under the assumptions of 
rationality (i.e. all modules are completely reducible) and $C_2$-cofiniteness. 
Our aim is to develop general theorems without the assumption of 
rationality so that we can use them 
in the aim of proving rationality of given VOAs. Our main target is an 
orbifold model. 
Throughout this paper, $T$ denotes a VOA with 
an automorphism $g$ of order $n$. Set 
$T^{(k)}=\{t\in T \mid  g(t)=e^{2\pi ik/n}t\}$, then 
$T=\oplus_{k=0}^{n-1}T^{(k)}$ and $T^g=T^{(0)}$ is a subVOA. 

Let $V=\oplus_{n=0}^{\infty}V_n$ be a VOA. 
In an orbifold module, we will study $T^{(0)}$ as a targeted VOA $V$. 
In \cite{M3} and \cite{M4}, the author have shown several results under 
the assumption that $V$ is $C_2$-cofinite. We will study some of them under weaker conditions.  

Let ${\rm Irr}(V)$ and ${\rm mod}(V)$ denote the set of irreducible $V$-modules and 
that of $\N$-graded $V$-modules with composition series of finite length.  
In this paper, for $U,W\in {\rm mod}(V)$, 
we will define a fusion product $U\boxtimes W$ as a projective limit in \S 3. 
As we will explain, $U\boxtimes W$ may not be a $V$-module, but  
$W\boxtimes U$ satisfies all conditions to be a $V$-module except for the  
lower truncation property and 
it is possible to define the $n$-th action $v_n$ of $v\in V$ on $W\boxtimes U$ 
and we are able to consider $V$-homomorphisms from $W\boxtimes U$ to $V$-modules. 

For an orbifold model, we will show: \\

\noindent
{\bf Proposition \ref{OSC}}\quad {\it  
Let $T$ be a VOA and $g\in {\rm Aut}(V)$ of order $n$. 
If all $T^{(k)}$ are $C_1$-cofinite as $V$-modules and $T^{(k)}\boxtimes T^{(1)}$ are  
$V$-modules, then all $T^{(k)}$ are simple currents.}  \\

Our definition of ``$C_1$-cofinite as a $V$-module" is slightly different from 
an ordinary definition of $C_1$-cofiniteness, see \S 2.3. 
A surjection $\alpha:D\to C$ of $V$-modules is called ``covering" of $C$ when 
$D$ is the unique submodule of $D$ satisfying $\alpha(D)=C$. We sometimes 
call $D$ a cover of $C$. The flatness property (i.e. preserving the exactness) of the tensor 
products is important in an ordinary ring theory. Unfortunately, in 
VOA theory, we need some assumption to expect it.  
We first pay attention to the following short exact sequences. \\

\noindent
{\bf Theorem \ref{FPP}}\quad {\it 
Let $V$ be a VOA and $\rho:P\to V$ a covering of $V$ as $V$-modules. 
If the weights of $\Ker(\rho)\boxtimes U$ are bounded below for $U\in {\rm mod}(V)$, then 
$$0\rightarrow {\rm \Ker}(\rho)\boxtimes U 
\xrightarrow{\epsilon \boxtimes {\rm id}_U} 
P\boxtimes U 
\xrightarrow{\rho\boxtimes {\rm id}_U} V
\boxtimes U\rightarrow 0  $$
is exact, where $\epsilon:\Ker(\rho)\to P$ is an embedding and ${\rm id}_U$ denotes 
the identity map on $U$.}\\

Since we do not assume rationality, we have to prepare for a non-semisimple category. 
One important concept we will introduce is ``semi-rigidity", which is a non-semisimple 
version of the rigidity. 

\begin{dfn}\label{SRD}
$\mbox{}\quad$ $W\in {\rm mod}(V)$ is called to be {\bf semi-rigid} 
if there is $\widetilde{W}\in {\rm mod}(V)$, $e_W\in \Hom_V(W\boxtimes \widetilde{W},V)$ and   
$e_{\widetilde{W}}\in \Hom_V(\widetilde{W}\boxtimes W,V)$ such that \\
(1) $W, \widetilde{W}$ are $C_1$-cofinite as $V$-modules and 
$W\boxtimes \widetilde{W}$ and $W\boxtimes (\widetilde{W}\boxtimes W)$ 
are $V$-modules, \\ 
(2) there are a $V$-module $Q$ and an embedding $\epsilon:Q\to \widetilde{W}\boxtimes W$ such that 
$e_{\widetilde{W}}\epsilon:Q\to V$ is a covering and \\
(3) 
in the diagram
$$\begin{array}{ccccc}
W\boxtimes Q&\xrightarrow{{\rm id}_W\boxtimes {\epsilon}}
&W\boxtimes (\widetilde{W}\boxtimes W)& \xrightarrow{{\rm id}_W\boxtimes e_{\widetilde{W}} } 
&W\boxtimes V \cr
&&\mbox{}\qquad\downarrow \mu  & & \cr
 &&(W\boxtimes \widetilde{W})\boxtimes W   & \xrightarrow{e_W\boxtimes {\rm id}_W} & V\boxtimes W
\end{array}\eqno{(1.1)}$$ 
$(e_W\boxtimes {\rm id}_W)\mu({\rm id}_W\boxtimes \epsilon)$ is surjective,  
where $\mu:W\boxtimes (\widetilde{W}\boxtimes W) \rightarrow (W\boxtimes \widetilde{W})\boxtimes W$ 
is a natural isomorphism for the associativity of 
products of intertwining operators (see (5.3)). 
\end{dfn}

Clearly, every simple current is semi-rigid, in particular, so is $V$. 
We will prove the following theorems: \\ 

\noindent
{\bf Theorem \ref{PDS}}\quad{\it 
If $C$ is a semi-rigid irreducible module 
and $D\xrightarrow{\alpha} C$ is a covering of $C$, 
then there is a covering $R\xrightarrow{\rho} V$ such that 
$D$ is a homomorphic image of $R\boxtimes C$. In particular, 
if $V$ has a projective cover $P_V$ and $P_V\boxtimes C$ is 
finitely generated, then $C$ has a projective cover $P_C$ which is 
isomorphic to a direct summand of $P_V\boxtimes C$. }\\

Therefore, if $V$ is projective as a $V$-module, then 
all semi-rigid irreducible modules are projective (Corollary 15). \\

\noindent
{\bf Theorem \ref{FLT}}\quad  
{\it Let $V$ be a $C_2$-cofinite VOA of CFT-type and $C$ a semi-rigid $V$-module. Then 
$$0\rightarrow {\rm \Ker}(\rho_C)\boxtimes W 
\xrightarrow{\epsilon\boxtimes {\rm id}_W} 
P_C\boxtimes W 
\xrightarrow{\rho_C\boxtimes {\rm id}_W} C
\boxtimes W\rightarrow 0  $$
is exact, where $0\rightarrow {\rm Ker}(\rho_C)\xrightarrow{\epsilon} P_C
\xrightarrow{\rho_C} C\to 0$ is a projective cover of $C$.  }\\

As corollaries, we have the following for a $C_2$-cofinite VOA $V$. \\
(1). If all $V$-modules are semi-rigid, then all $V$-modules are flat, that is, 
fusion products preserve the exactness.  \\
(2). If all $V$-modules are semi-rigid and $V$ is projective as a $V$-module, 
then $V$ is rational. [Corollary \ref{AP}]\\  

Since $V$ has integer weights, the investigations of its extensions are easier 
than those of other modules. For example, we will show the following for any VOA $V$. \\
(3).  If a simple VOA $V$ contains a rational subVOA with the same Virasoro element, 
then $V$ is projective as a $V$-module. [Corollary \ref{SR}]\\

We first introduce a weaker condition than $V'\cong V$. \\

\noindent
Condition I:  \quad {\it Let $V$ be a simple VOA and 
for each $W\in {\rm mod}(V)$, there is $\widetilde{W}\in {\rm mod}(V)$ such that 
$\Hom_V(W\boxtimes \widetilde{W},V)\not=0$.} \\

If $V'\cong V$, then we can take a restricted dual $W'$ of $W$ as $\widetilde{W}$. 
We remark that if products of intertwining operators satisfies the associativity, then 
since there is an epimorphism $W\boxtimes W'\to V'$ by adjoint operators and 
$(W\boxtimes W')\boxtimes Q\cong W\boxtimes (W'\boxtimes Q)$, 
Condition I comes from the following:\\

\noindent
Condition I':\quad {\it There is $Q\in {\rm mod}(V)$ such that $\Hom_V(V'\boxtimes Q,V)\not=0$.}\\

Our next aim is to find a condition under which all $V$-modules become semi-rigid. 
Let us consider a trace function of $v$ on an irreducible $V$-module 
$W=\oplus_{n=0}^{\infty} W_{n+r}$: 
$$\Psi_{W}(v,\tau)=\Tr_W (o(v) q^{(L(0)-c/24)})
=\sum_{n=0}^{\infty}(\Tr_{W_{r+n}}o(v))q^{(n+r-c/24)} \qquad \mbox{ for }v\in V$$
and define a map 
$$\Psi_{W}(\ast,\tau): v\in V \mapsto  \Psi_W(v,\tau),$$
where $q$ denotes $e^{2\pi i\tau}$, $c$ is a central charge of $V$ and 
$o(v)$ denotes a grade-preserving operator of $v$ on $W$, 
(e.g. $o(v)=v^M_{\wt(v)-1}$ for $v\in V_{\wt(v)}$ and $Y^M(v,z)=\sum v_n^Mz^{-n-1}$). 
Abusing the notation, we will also call $\Psi_W$ trace functions. 
We then consider its $S$-transformation $S(\Psi_{V})$ (corresponding to 
$\begin{pmatrix}0&-1\cr 1&0\end{pmatrix}\in SL_2(\Z)$). 
It satisfies  
$$S(\Psi_{V})(v,\tau)=(\frac{-1}{\tau})^{\wt(v)}\Psi_V(v,-1/\tau)$$
for $v\in V_{\wt(v)}$ with $L(1)v=L(2)v=0$. 
As the author has shown in \cite{M1}, if $V$ is $C_2$-cofinite, then $S(\Psi_{V})$ equals 
a linear combination of trace functions and pseudo-trace functions. 
In the case where $V$ is rational, $S(\Psi_{V})$ has no pseudo-trace functions, 
as Zhu has shown in \cite{Zh}.

We will consider the following condition: \\

\noindent
Condition II: \quad $S(\Psi_V)$ is a linear combination of 
trace functions on $V$-modules. \\

Condition II looks strong, but for an orbifold model, we have: \\

\noindent
{\bf Theorem \ref{NST}} \quad {\it 
Let $T$ be a VOA of CFT-type with a finite automorphism $g$ and assume that 
a fixed point subVOA $T^g$ is $C_2$-cofinite. 
If $T$ satisfies Condition II, so does $T^g$. }\\

The key result for orbifold model is the following: \\

\noindent
{\bf Theorem \ref{NCE}} \quad {\it
If $V$ satisfies Condition I and Condition II, 
then all simple $V$-modules are semi-rigid. 
Furthermore, if $V$ is a rational VOA of CFT-type satisfying Condition I, 
then $\Psi_U$ has nonzero coefficient in $S(\Psi_V)$ for every simple $V$-module $U$.}\\

As a corollary, we will prove: \\

\noindent
{\bf Corollary \ref{OBF}}\quad {\it 
Let $g$ be a finite automorphism of a VOA $T$ and assume that 
a fixed point subVOA $V:=T^g$ is a $C_2$-cofinite of CFT-type and satisfies Condition I, 
If $T$ is rational and $V$ satisfies Condition I, 
then $V$ is rational. Moreover, for a simple $V$-module $W$, there are $k\in \Z$ and 
a $g^k$-twisted $T$-module $U$ such that $W$ is a $V$-submodule of $U$. 
}\\

We organize this paper in the following way. 
We treat a fusion produce as a projective limit of 
intertwining operators in Section 3. In Section 5, as an application of fusion products, 
we present a product of two intertwining operators by a homomorphism. This will be one of 
the main tools in this paper. 
In Subsection 4.1 we show that the principal projective cover preserves the exactness with respect to 
the fusion products. In Section 6, we study an orbifold model under weaker conditions than $C_2$-cofiniteness. 
We study three transformations in Section 8 and prove that all modules are semi-rigid under some condition. 
In Section 9 and 10, we prove an orbifold conjecture under the assumption that 
a fixed point subVOA is $C_2$-cofinite. \\

\noindent
{\bf Acknowledgement} \\
The author would like to express special thanks to H.~Yamauchi and T.~Arakawa 
for giving me the chance to explain the detail of his proofs.  
He thanks T.~Abe, K.~Harada, A.~Matsuo, K.~Naoi, H.~Yamada and K.~Tanabe 
for their pointed questions.   
He also thanks to G.~H\"{o}hn, V.~Nikulin and N.~Scheithauer, the organizers of 
the conference held at Edinburgh in 2009 where he was stimulated to extend 
his results to general cases. 

\section{Preliminary results}
\subsection{Notation}
Throughout this paper, $V$ denotes a VOA 
$(V=\oplus_{i=0}^{\infty},Y(\cdot,z),\1,\omega)$, 
where $\omega$ is a Virasoro element, 
$\1$ is the vacuum, and $Y(v,z)=\sum_{n\in \Z} v_nz^{-n-1}\in 
\End(V)[[z,z^{-1}]]$ denotes a 
vertex operator of $v\in V$ satisfying the conditions 
$1\sim 4$ in Definition \ref{LOG} by replacing all $W,U,R$ by $V$. 
We also have $Y(\1,z)={\rm id}_V$, $\lim_{z\to 0}Y(v,z)\1=v$ and the coefficients of $Y(\omega,z)$ 
satisfy the Virasoro algebra relations. If $W=V$, $U=R=M$ in Definition \ref{LOG} and a set  
$$\{Y^M(v,z)=\sum_{n\in \Z} v^M_nz^{-n-1}\in 
\End(M)[[z,z^{-1}]] \mid v\in V\}$$ 
satisfies the same conditions $1\sim 4$, then $M$ is called 
a $V$-module (including weak modules). 
We use notation $``\wt"$ (weight) to denote the semisimple part of $L(0)$ and 
$L(0)^{nil}$ to denote $L(0)-{\wt}$. 
All modules in this paper are $\C$-graded $M=\oplus_{r\in \C} M_r$ by weights 
excepted as otherwise noted, that is, 
they are direct sums of generalized eigenspaces $M_r$ of $L(0)$ with eigenvalue $r$. 
For a $\C$-graded module $M=\oplus_{n\in \C}M_r$, 
$M'$ denotes the restricted dual $V$-module 
$\oplus_{r\in \C}\Hom(M_r,\C)$, where an adjoint vertex operator 
$Y^{M'}(v,z)$ on $M'$ is given by 
$$\langle Y^{M'}(v,z)w',w\rangle
=\langle w',Y(e^{zL(1)}(-z^{-2})^{L(0)}v,z^{-1})w\rangle $$ 
for $w'\in M'$ 
and $w\in M$, where $\langle w',w\rangle$ denotes $w'(w)$, see \cite{FHL}.

\subsection{(logarithmic) intertwining operators}
Similar to a vertex operator 
$Y^M(v,z)=\sum v^M_nz^{-n-1}\in \End(M)[[z,z^{-1}]]$ on a module $M$, 
it is natural to consider an intertwining operator from 
a module $U$ to another module $R$ as a formal power series (\cite{TK}). 
However, without the assumption of rationality, 
there is no reason for all intertwining operators to have specific 
forms like formal power series. As Huang has shown in \cite{HD}, if 
a $V$-module $W$ is $C_1$-cofinite, then all intertwining operators 
$\CY$ of type $\binom{R}{W\,\, U}$ 
has a shape of formal power series with $\log z$ terms: 
$$\CY(w,z)=\sum_{i=0}^K\sum_{m\in {\mathbb C}} 
w^{\CY}_{(i,m)}z^{-m-1}\log^i z\in \Hom(U,R)\{z\}[\log z] $$ 
and is called ``logarithmic type" (see \cite{Mil},\cite{G} and \cite{F}). 
As long as we have a possibility to treat non-semisimple modules, 
it is necessary to consider all of such operators. 
In this paper, we will call them intertwining operators, too.

\noindent
\begin{dfn}\label{LOG} Let $W$, $U$ and $R$ be 
 $\N$-graded $V$-modules. 
A (logarithmic) intertwining operator of type 
$\binom{R}{W\,\,U}$ is a linear map 
$$\begin{array}{l}
  \CY(,z):W \rightarrow \Hom(U,R)\{z\}[\log(z)] \cr
  \displaystyle{\CY(w,z)=\sum_{i=0}^k\sum_{m\in \C} 
w_{(i,m)}z^{-m-1}\log^iz} 
\end{array}$$
satisfying the following conditions: \\
1. The lower truncation property: for each $u\in U$ and $i$,  
$w_{(i,r)}u=0$  for ${\rm Re}(r)>\!\!>0$, 
where ${\rm Re}(r)$ denotes the real part of $r\in \C$.  \\
 2. $L(-1)$-derivative property: 
$\CY(L(-1)w,z)=\frac{d}{dz}\CY(w,z)$ for $w\in W$. \\
3. Commutativity: $v_n^R\CY(w,z)-\CY(w,z)v_n^U=
\displaystyle{\sum_{i=0}^{\infty}} \binom{n}{i}\CY(v_i^Ww,z)z^{n-i} 
\mbox{  for  }v\in V$. \\
4. Associativity: for $v\in V$ and $n\in \N$, \\
$\CY(v^W_{n}w,z)=
\displaystyle{\sum_{m=0}^{\infty}}(-1)^m\binom{n}{m}v^R_{-m-1}z^m\CY(w,z)
+\CY(w,z)
\displaystyle{\sum_{m=0}^{\infty}}(-1)^{m+n+1}\binom{n}{m}v^U_mz^{-m-1}.$ 
\end{dfn}

\noindent
Here and henceforth, 
$Y^X(v,z)=\sum v^X_nz^{-n-1}$ denotes a vertex operator of $v\in V$ on 
a $V$-module $X$ and $\CI\binom{R}{W\,\,U}$ denotes 
the set of intertwining operators of type 
$\binom{R}{W\quad U}$. 

Let $\CY^{(i)}(w,z):=\sum_{n\in \C} w_{(i,n)}z^{-n-1}$ denote the coefficient of 
$\log^iz$ in $\CY(w,z)$ for $i=0,1,...,K$, 
If $W$ is finitely generated, then there is an integer $K$ such that 
$\CY^{(n)}(w,z)=0$ for any $n>K$ and $w\in W$. 
Moreover, since vertex operators $Y^M(v,z)$ on modules $M$ 
have no $``{\log} z"$ terms, 
every $\CY^{(i)}$ satisfies all properties of intertwining 
operators except the $L(-1)$-derivative property. 
From the $L(-1)$-derivative property for $\CY$, 
we have two important properties:
$$\begin{array}{rl}
\displaystyle{\CY^{(m)}(w,z)=}&\!\!\displaystyle{\frac{1}{m!}(z\frac{d}{dz}-zL(-1))^m\CY^{(0)}(w,z)}, 
\qquad \mbox{ and }\cr
\displaystyle{(i+1)w_{(i+1,n)}u
=}&\!\!\displaystyle{-L(0)^{nil}w_{(i,n)}u+(L(0)^{nil}w)_{(i,n)}u
+w_{(i,n)}(L(0)^{nil}u)}.  
\end{array} \eqno{(2.1)}$$
In particular, we have 
$$(z\frac{d}{dz}-zL(-1))^{K+1}\CY^{(i)}(w,z)=0 \quad \mbox{ for any }i. \eqno{(2.2)}$$
On the other hand, 
for such a formal power series $\CY^0$ satisfying all conditions 
in Definition \ref{LOG} except $L(-1)$-derivative property but (2.2),
$$\CY(w,z)={\displaystyle \sum_{i=0}^{K}} \left\{\frac{1}{i!}(zL(-1)-z\frac{d}{dz})^i\CY^0(w,z)\right\}\log^iz 
\eqno{(2.3)}$$
is an intertwining operator. We note that 
$\sum_{i=0}^K\CY^{(i)}(w,ze^{2\pi i})(\log z+1)^i$ 
is also an intertwining operator for 
$\sum_{i=0}^K\CY^{(i)}(w,z)\log^iz\in \CI\binom{\ast}{W\,\,\ast}$.

We also note that if $L(0)$ acts on $W$ semisimply and $\CY\in \CI\binom{U}{W\,\,U}$, then 
since $L(0)^{nil}$ is a nilpotent operator and commutes with all grade-preserving operators, 
the trace of $w_{(i,\wt(w)-1)}q^{L(0)}$ on $U$ is zero for $i\geq 1$ by (2.1), 
where $\CY(w,z)=\sum_{i,n} w_{(i,n)}z^{-n-1}\log^iz$. 
Therefore, the followings come from (2.1). \\

\begin{lmm} 
(1) If $L(0)$ acts on $W$ semisimply, then 
$\Tr_U\CY(w,z)q^{L(0)}=\Tr_U\CY^{(0)}(w,z)q^{L(0)}$ for $\CY\in \CI_{W,U}^U$. \\
(2) If $L(0)$ acts on $W,U,R$ semi-simply, then 
every $\CY\in \CI\binom{R}{W\,\,U}$ is a formal power series. \end{lmm}

\subsection{$C_m$-cofiniteness for a module}
For a $V$-module $W$ and $m\in \N$, set 
$$C_m(W)=\langle v_{-m}w\mid v\in V, \wt(v)>0, w\in W\rangle. $$
If $\dim W/C_m(W)<\infty$, then we call $W$ to be $C_m$-cofinite as a $V$-module. 
Our $C_1$-cofiniteness is slightly different from the ordinary definition of $C_1$-cofinite. 
For example, $V$ is always $C_1$-cofinite as a $V$-module. 
Among these finiteness conditions, $C_2$-cofiniteness is the most important and
offers many nice properties. 
For example, we have:

\begin{prn}\label{C2Pr}
Let $V$ be a $C_2$-cofinite VOA. Then we have the followings: \\
(i) Every weak module is $\Z_+$-gradable and weights are all rational numbers {\rm \cite{M2}}. \\
(ii) Evey $n$-th Zhu algebra $A_n(V)$ is finite dimension and
the number of inequivalent simple modules is finite, 
{\rm \cite{GN},\cite{DLM}}.\\
(iii) Set $V=B+C_2(V)$ for 
a finite dimensional subspace $B$ spanned by homogeneous elements. 
Then for any weak module $W$ generated from one element $w$ has the 
following spanning set $\{v^1_{n_1}....v^k_{n_k}w 
\mid v^i\in B, 
\quad  n_1<\cdots <n_k\}$. In particular, 
every $V$-module is $C_n$-cofinite as a $V$-module for 
any $n=1,2,...$, {\rm \cite{M2},\cite{Bu},\cite{GN}}. 
\end{prn}

We note that except Condition I, 
all conditions in this paper are satisfied by a simple $C_2$-cofinite VOA and 
its simple modules.

\section{Fusion products}
In this section, for $W, U\in {\rm mod}(V)$ we explain  
a fusion product $W\boxtimes U$ defined by 
intertwining operators as a projective limit. 
We call an intertwining operator 
$$\CY(w,z)=\sum_{i=0}^K\sum_{m\in \C} w_{(i,m)}z^{-m-1}\log^iz$$ 
of type $\binom{T}{W\quad U}$ ``surjective" 
if the images of coefficients of $\CY$ spans $T$, that is, 
$$\langle w_{(i,m)}u\mid w\in W, u\in U, m\in \C, i=0,...,K\rangle_{\C}=T.$$ 

Consider the set of surjective intertwining operators of $W$: 
$$ \CF(W,U)=\{(\CY,F)\mid  F\in {\rm mod}(V), \CY\in \CI\binom{F}{W\,U}
\mbox{ is surjective}\} \eqno{(3.1)}$$
and introduce an equivalent relation 
$(\CY_1,F^1)\cong (\CY_2,F^2)$ if there is an 
isomorphism $\sigma:F^1\to F^2$ such that $\CY_2=\sigma\CY_1$. 
We define a partial order $\leq $ in $\CF(W,U)/\cong$ as follows: \\
For $(\CY_1,F^1),(\CY_2,F^2)\in \CF(W,U)$, \\ 
$\CY_1\leq \CY_2$ $\Leftrightarrow$ ${}^{\exists}f\in \Hom_V(F^2,F^1)$ s.t. $f\CY_2=\CY_1$. 

Clearly, if $\CY_1\leq \CY_2$ and $\CY_2\leq \CY_1$, then we 
have $(\CY_1,F^1)\cong (\CY_2,F^2)$.   

\begin{lmm}\label{DS}
$\CF(W,U)/\cong$ is a (right) directed set. 
\end{lmm}

\pr 
\hspace{-4mm} For any $\CY_1,\CY_2 \in \CF(W,U)$, say 
$\CY_1\in \CI\binom{F^1}{W\,\,U}$ and $\CY_2\in \CI\binom{F^2}{W,U}$, 
we define $\CY$ by 
$$\CY(w,z)u=(\CY_1(w,z)u, \CY_2(w,z)u)\in (F^1\times F^2)\{z\}[\log z] \quad 
\mbox{ for }w\in W, u\in V.$$
Clearly, $\CY$ is an intertwining 
operator of type $\binom{F^1\times F^2}{W\,\,U}$. 
Let $F\subseteq F^1\times F^2$ denote the subspace 
spanned by all images of coefficients of $\CY(w,z)u$ with $w\in W, u\in U$, 
then since $F_1$ and $F_2$ have composition series of finite length, 
so does $F$ and so $(\CY,F)\in \CF(W,U)$. 
By the projections $\pi_i:F^1\times F^2 \rightarrow F^i$, we have 
$\pi_1(\CY)=\CY_1$ and $\pi_2(\CY)=\CY_2$, that is, 
$(\CY_1,F^1)\leq (\CY,F)$ and $(\CY_2,F^2)\leq (\CY,F)$ as we 
desired. 
\prend

Set $\CF(W,U)=\{(\CY_i,F^i) \mid i\in I\}$ and consider the product 
$(\prod_i \CY_i, \prod_i F^i)$. Let $F_{(r)}$ be a subspace of $\prod T_i$ spanned 
by all coefficients 
$\prod w^{\CY_i}_{\wt(w)-1-r+\wt(u)}u$ of weights $r$ for homogeneous 
elements $w\in W$ and $u\in U$ and we set 
$F=\coprod_{r\in \C} F_{(r)}$ and we consider $(F,\prod_i \CY^i)$. 

\begin{dfn}\label{FPT} 
We call $F$ a ``fusion product" of $W$ and $U$ and we denote it by 
$W\boxtimes U$. 
\end{dfn}

Namely, we define a fusion product as a projective limit of $\CF(W,U)/\cong$. 
We note that $U\boxtimes W$ may not be a $V$-module, but  
$W\boxtimes U$ satisfies all conditions to be a $V$-modules except for the  
lower truncation property. 
Hence if the set of conformal weights of $F^i$ in $\CF(W,U)$ has a lower bound, 
then $W\boxtimes U$ is a $V$-module and 
$\prod \CY^i$ is a surjective intertwining 
operator of type $\binom{W\boxtimes U}{W\quad U}$.  
We also note that even if $W\boxtimes U$ is not a $V$-module, 
it is possible to define the actions $v_n$ on $W\boxtimes U$ and we are able to 
consider $V$-homomorphisms from 
$W\boxtimes U$ to $V$-modules. 
In any case, $U\boxtimes W$ is $\C$-graded by the construction.  
Throughout this paper, we fix one surjective intertwining operator 
$\CY_{A,B}^{\boxtimes}\in \CI\binom{A\boxtimes B}{A\,\,B}$ 
for each pair $A,B\in {\rm mod}(V)$.

We next show:

\begin{prn}\label{EOF} Let $W,U\in {\rm mod}(V)$. If an $n$-th Zhu bi-module $A_n(W)$ is of finite dimensional 
for any $n\in \N$, then 
$\CF(W,U)$ contains a unique maximal element. 
In other words, a fusion product $W\boxtimes U$ is well-defined as an element of ${\rm mod}(V)$. 
\end{prn}

Before we start the proof, we give a brief review of an $n$-th 
Zhu algebra $A_n(V)$ and an $n$-th Zhu bi-module $A_n(W)$ from a view point of operators. 
For $\N$-graded modules $U=\oplus_{i=0}^{\infty} U_{i+r}$ and 
$F=\oplus_{i=0}^{\infty}F_{i+t}$ with lowest weights $r$ and $t$, respectively, 
we restrict our interest into operators 
$$o_{t-r}^{\CY}(w)=w^{\CY}_{(0,\wt(w)-1+r-t)}:\oplus_{j=0}^nU_{j+r}\to 
\oplus_{j=0}^nF_{j+t},$$ 
where $\CY(w,z)=\sum_{i=0}^K\sum_{m\in \C}w^{\CY}_{(i,m)}z^{-m-1}\log^iz\in 
\CI\binom{F}{W\,\,U}$ of $w\in W$ and $o^{\CY}_s(w)$ denotes an operator of $w$ 
in $\CY^{(0)}(w,z)$ shifting grade by $s$.  
From Associativity and Commutativity, 
for $v\in V$ and $w\in W$, we can find $v\ast w$ and $w\ast v$ in $W$ such that 
$o^{\CY}_{t-r}(v\ast w)=o_0(v)o^{\CY}_{t-r}(w)$ and 
$o^{\CY}_{t-r}(w\ast v)=o^{\CY}_{t-r}(w)o_0(v)$ for any 
$\CY\in \CI\binom{F}{W\,\,U}$, where $o_0(v)$ denotes a grade preserving operator of $v$.  If we set 
$$O_n(W)=\{w\in W\mid o^{\CY}_{t-r}(w)=0 \mbox{ for any }
\CY\in \CI\binom{F}{W\,\,U}, 
U, F\in {\rm Ind}(V) \}$$
where ${\rm Ind}(V)$ is the set of indecomposable $V$-modules in ${\rm mod}(V)$ and 
$r$ and $t$ are the lowest weights of $U$ and $F$, respectively, 
then $A_n(V)=V/O_n(V)$ is an associative algebra and $A_n(W)=W/O_n(W)$ becomes 
an $A_n(V)$-bi-module.  \\

\noindent
{\bf [Proof of Proposition \ref{EOF}]} \qquad 
We may assume that $W$ and $U$ are indecomposable. 
Since the actions of $V$ shift the weights by integers, there are 
$r,s\in \C$ such that $W=\oplus_{n=0}^{\infty}W_{n+s}$ and 
$U=\oplus_{n=0}^{\infty}U_{n+r}$. Let $(\CY,F)\in \CF(W,U)$ 
with $\CY(w,z)=\sum_{i=0}^K\sum_m w_{(i,m)}z^{-m-1}\log^iz$.  
We will show that the lengths of composition series of $F$ have an upper bound.    
Since $\dim A_0(V)<\infty$, there are only finitely many conformal weights 
and so we assume $F=\oplus_{n=0}^{\infty} F_{n+t}$ with a lowest weight $t$. 
Moreover, for a sufficiently large $N$, 
a composition series of $F$ as a $V$-module corresponds to a 
composition series of $F_{N+t}$ as an $A_N(V)$-module.  
Since $\CY$ is surjective, we have 
$$\langle w_{(i,\wt(w)-1-N+n+r-t)}(U_{n+r})\mid w\in W, n\in \N,i=0,\ldots,K\rangle_{\C}
=F_{N+t}.$$ 
As we did in the explanation of an $N$-th Zhu algebra, we have 
$$\langle w_{(i,\wt(w)-1+r-t)}(U_{N+r})\mid w\in W, i=0,\ldots,K\rangle_{\C}
=F_{N+t}.$$ 
We have to note that for $w_{(i,m)}$ with $i\not=0$, we are able to choose another $\widetilde{\CY}$ 
such that $\widetilde{\CY}^{(0)}=\CY^{(i)}$ as we mentioned at (2.3). 
Therefore, there is a surjection  
$$\phi:\coprod_{h=0}^K\left(A_N(W)\otimes U_{N+r}\right) \rightarrow F_{N+r}$$ 
given by $\phi(\coprod_{h=0}^K(w^{(h)}\otimes u^{(h)}))
=\sum_{h=0}^K w^{(h)}_{(h,\wt(w^{(h)})-1+r-t)}u^{(h)}$ and so we have 
$$\dim T_{N+r}\leq \dim A_N(W)\times \dim U_{N+r}\times (K+1). $$
The right hand side does not depend on the choice of $F$. 
Therefore, $\CF(W,U)/(\cong)$ has a unique maximal element, which 
is a fusion product.   
\prend

We next explain about induced homomorphisms among fusion products. 
For a homomorphism $\phi: A\rightarrow B$ of $V$-modules, a formal operator 
$\CY\in \CI\binom{B\boxtimes D}{A\,\,D}$ defined by 
$$\CY(a,z)d=\CY_{B,D}^{\boxtimes}(\phi(a),z)d$$ 
becomes an intertwining 
operator of type $\left({}_{A\quad D}^{B\boxtimes D}\right)$. Therefore, by the 
maximality of fusion products, there is a $V$-homomorphism 
$\phi\boxtimes {\rm id}_D: A\boxtimes D \rightarrow B\boxtimes D$ such that 
$$\phi\boxtimes {\rm id}_D\cdot(\CY(a,z)d)=\CY(\phi(a),z)d.$$  
Similarly, we can define ${\rm id}_D\boxtimes 
\phi:D\boxtimes A\rightarrow D\boxtimes B$. The following shows that 
fusion products preserve the right exactness of sequences.  

\begin{prn}\label{RightFl}  Let $A, B, C, D\in {\rm mod}(V)$. 
If 
$$A\xrightarrow{\phi} B\xrightarrow{\sigma} C \rightarrow 0$$ 
is exact, then so is 
$$ A\boxtimes D \xrightarrow{\phi\boxtimes {\rm id}_D} 
B\boxtimes D 
\xrightarrow{\sigma\boxtimes {\rm id}_D} C\boxtimes D \rightarrow 0.$$
\end{prn}

\pr
Clearly, $(\sigma\boxtimes {\rm id}_D)\cdot(\phi\boxtimes {\rm id}_D)
=(\sigma\cdot \phi)\boxtimes {\rm id}_D=0$. Since we may view 
$(\CY,F)\in \CF(C,D)$ as $(\CY(\sigma),F)\in \CF(B,D)$, 
$\sigma\boxtimes {\rm id}_D$ is surjective and so we may view 
$(C\boxtimes D)'\subseteq (B\boxtimes D)'$. 
We may also assume $A\subseteq B$ and $C=B/A$. Consider a canonical bilinear pairing 
$$\langle w , \CY_{B,D}^{\boxtimes}(b,z)d\rangle \in \C\{z\}[\log z] $$
for $w\in (B\boxtimes D)', b\in B$ and $d\in D$. 
Clearly, if $w\in (C\boxtimes D)'$, then 
$$\langle w, \CY_{B,D}^{\boxtimes}(a,z)d\rangle =0 \eqno{(3.2)}$$
for any $a\in A$. 
On the other hand, if $w\in (B\boxtimes D)'$ satisfies (3.2) 
for all $a\in A$, then $\langle w, \CY_{B,D}^{\boxtimes}(b,z)d\rangle $ is well 
defined for $b\in B/A=C$. Therefore, 
$$0 \rightarrow (C\boxtimes D)' \rightarrow (B\boxtimes D)' 
\rightarrow (A\boxtimes D)'$$
is exact, which implies that $\boxtimes D$ preserves the right exactness. 
\prend \vspace{-4mm}

\section{Projective covers}
Let us start with an explanation of projective modules and 
projective covers. 

\begin{dfn} A $V$-module $P\in {\rm mod}(V)$ is called ``projective" 
if every $V$-epimorphism $f:W\rightarrow P$ splits for $W\in {\rm mod}(V)$. 
If $P\to D$ is a cover and $P$ is projective, we call $P$ a ``projective cover" of $D$. 
\end{dfn}

Different from an ordinary ring theory, 
$V$ is not necessarily projective as a $V$-module. 

\begin{prn}\label{Covering}  Let $V$ be a VOA and $D, U\in {\rm mod}(V)$. 
Then we have: \\
(1) If $D\in {\rm mod}(V)$ is not projective, then 
there is a proper covering $F\xrightarrow{\beta}D$. \\
(2) If $\dim A_n(V)<\infty$ for any $n$, then every $U\in {\rm mod}(V)$ has 
a projective cover. 
\end{prn}

\pr 
Proof of (1). Since $D$ is not projective, 
there is a non-split surjection $\beta:F\rightarrow D$. We choose $F$ minimal 
among such non-split extensions. Clearly, $F$ is a covering. \\
Proof of (2). 
Since $U$ is finitely generated, there are $u^{(1)},\ldots,u^{(k)}$ in $U$ such that 
$U=Vu^{(1)}+\cdots+Vu^{(k)}$, where $Vu=\{v_mu\mid v\in V, m\in \Z\}$ 
denotes a submodule generated from $u$ by \cite{DM}. 
We may assume that $u^{(i)}$ are all homogeneous with the same weight, say $r$.  
Let $D\xrightarrow{\alpha}U$ be a covering of $U$. 
Since $\dim A_0(V)<\infty$, the number of irreducible $V$-modules is finite and 
so there is an integer $N$ which does not depend on $D$ such that 
the conformal weight of $D$ is greater than $r-N$. 
Choose homogeneous elements $d^{(i)}\in D$ so that $\alpha(d^{(i)})=u^{(i)}$ for 
every $i$. Since $\alpha$ is a covering, $D=Vd^{(1)}+\cdots+Vd^{(k)}$. 
Then the subspace $D_r$ of $D$ of weight $r$ 
is spanned by $\{o(v)d^{(i)}\mid v\in V, i=1,...,k\}$ and so 
$\dim D_r\leq k\dim A_N(V)$. Thus, the length of composition series of 
coverings of $U$ is bounded and so $U$ has a projective cover by (1). 
\prend

\subsection{Principal projective cover and fusion products}
In this subsection, we will prove the following:

\begin{thm}\label{FPP} 
Let $V$ be a VOA and $\rho:P\to V$ a covering of $V$. 
Then for $W\in {\rm mod}(V)$, 
$$0\rightarrow {\rm \Ker}(\rho)\boxtimes W 
\xrightarrow{\epsilon \boxtimes {\rm id}_W} 
P\boxtimes W 
\xrightarrow{\rho\boxtimes {\rm id}_W} V
\boxtimes W\rightarrow 0  \eqno{(4.1)}$$
is exact if the weights of $\Ker(\rho)\boxtimes W$ are bounded below, 
where $\epsilon:\Ker(\rho)\to P$ is an embedding. 
\end{thm}

Before we start the proof, we will prove the following lemma:

\begin{lmm}\label{KOP}
Let $\rho:P\to V$ be a covering of $V$ and 
$p\in P$ a homogeneous element satisfying $\rho(p)=\1$. 
Then 
$$\Ker(\rho)=\langle v_np\in P \mid v\in V, n\geq 0\rangle_{\C}. \eqno{(4.2)}$$
\end{lmm}

\pr
Since 
$$Q=\langle v_np \mid v\in V, n\in \Z \rangle_{\C}$$
is a submodule of $P$ by \cite{DM} and $\rho(Q)=V$, 
we have $Q=P$. 
Furthermore, since $(L(-1)v)_mp=-mv_{m-1}p$, we have 
$$P=\langle v_np \mid v\in V, n\geq -1\rangle_{\C}
=\langle v_{-1}p\mid v\in V\rangle_{\C}
+\langle v_mp\mid v\in V, m\in \N\rangle_{\C}.$$ 
Since $\rho: \langle v_{-1}p\mid v\in V\rangle_{\C}\to 
\langle v_{-1}\1 \mid v\in V\rangle$ is an injection and 
$\rho(\langle v_mp\mid v\in V, m\in \N\rangle_{\C})=
\langle v_m\1\mid v\in V, m\in \N\rangle_{\C}=0$, 
we have the desired result. 
\prend

Let us start the proof of Theorem \ref{FPP}. \\

\noindent
$[${\bf  Proof of Theorem \ref{FPP}}$]$\quad 
Set $Q=\Ker(\rho)$ and $\epsilon:Q\to P$ denotes an embedding. 
By Proposition \ref{RightFl}, it is sufficient to show that 
$\epsilon\boxtimes {\rm id}_W:Q\boxtimes W\to P\boxtimes W$ is an 
injection.  
Let $\CY_{Q,W}^{\boxtimes}(a,z)=\sum_{j=0}^K\sum_m a_{(j,m)}z^{-m-1}\log^jz$ 
be a fusion product intertwining operator.  
Set $I=\CY^{(0)}$, that is, 
$I(a,z)=\sum_{r\in \C} a_{(0,r)}z^{-r-1}$ and we denote $a_{(0,r)}$ by $a^I_r$. 
We note that 
$I(\ast,z)$ is an $(z\frac{d}{dz}-L(-1)z)$-nilpotent intertwining operator as we 
showed at (2.2).  

We choose a homogeneous element $p\in P$ satisfying $\rho(p)=\1$. 
Consider a vector space 
$$R=p_{-1}W\oplus (Q\boxtimes W),$$ 
where $p_{-1}w$ is a formal element and $p_{-1}W=\{p_{-1}w \mid w\in W\}$. 
We view $Q\boxtimes W$ as a $V$-submodule of $R$ and define a $V$-module structure on $R$.  
Define an action $v^R_n$ of $v$ on $R=p_{-1}W\oplus Q\boxtimes W$ by 
$$v^R_n(p_{-1}w)
=p_{-1}(v_nw)+\sum_{i=0}^{\infty}\binom{n}{i}(v_ip)^I_{-1+n-i}w. $$ 
We note that since $v_ip\in Q$ for $i\geq 0$, $v_n^R$ is well-defined. 
By the direct calculations, we have $L(-1)$-Derivation:
$$\begin{array}{rl}
\displaystyle{(\omega_0v)^R_n(p_{-1}w)}=
&\displaystyle{p_{-1}((\omega_0v)_nw)+\sum_{i=0}^{\infty}
((\omega_0v)_ip)^I_{-1+n-i}w}\cr
=&\displaystyle{-np_{-1}(v_{n-1}w)+\sum_{i=0}^{\infty}
\binom{n}{i}(-iv_{i-1}p)^I_{-1+n-i}w}
=\displaystyle{-nv^R_{n-1}(p_{-1}w)}. 
\end{array} \eqno{(4.3)}$$
We also have Commutativity:
$$\begin{array}{l}
(v^R_mu^R_n-u^R_nv^R_m)(p_{-1}w) \cr
=\displaystyle{p_{-1}((v_mu_n-u_nv_m)w)+
\sum_{i=0}^{\infty}\binom{m}{i}(v_ip)^I_{m-1-i}\ast u_nw
-\sum_{i=0}^{\infty}\binom{n}{i}(v_ip)^I_{m-1-i}u_nw}\cr
\mbox{}\qquad+\displaystyle{\sum_{i=0}^{\infty}\binom{n}{i}(u_ip)^I_{n-1-i}v_mw
-\sum_{i=0}^{\infty}\binom{n}{i}(u_ip)^I_{n-1-i}v_mw }\cr
\mbox{}\qquad+\displaystyle{\sum_{i=0}^{\infty}\sum_{j=0}^{\infty}
\binom{n}{i}\binom{m}{j}(v_ju_ip)^I_{m+n-1-i-j}w
-\sum_{i=0}^{\infty}\sum_{j=0}^{\infty}\binom{m}{j}\binom{n}{i}
(u_iv_jp)^I_{m+n-1-i-j}w }\cr
=\displaystyle{p_{-1}([v_m,u_n]w)+\sum_{i=0}^{\infty}\sum_{j=0}^{\infty}
\binom{n}{i}\binom{m}{j}([v_i,u_j]p)^I_{m+n-1-i-j}w }\cr
=\displaystyle{p_{-1}([v_m,u_n]w)
+\sum_{i=0}^{\infty}\sum_{j=0}^{\infty}\sum_{k=0}^{\infty}
\binom{n}{i}\binom{m}{j}\binom{i}{k}((v_ku)_{i+j-k}p)^I_{m+n-1-i-j}w }\cr
=\displaystyle{p_{-1}(\sum_{j=0}^{\infty}\binom{m}{j}(v_ju)^I_{m+n-j}w)
+\sum_{j=0}^{\infty}\sum_{i=0}^{\infty}\binom{m}{j}
\binom{n+m-j}{i}((v_ju)_ip)^I_{m+n-1-j-i} }w\cr
=\displaystyle{\sum_{j=0}^{\infty}\binom{m}{j}(v_ju)^R_{m+n-j}(p_{-1}w)}. \hfill (4.4)
\end{array}$$
By solving $(v_nu)_r$ from Commutativity for $n\geq 0$, we have Associativity: 
$$ \sum_{i=0}^{\infty}\binom{n}{i}(-1)^i\{v^R_{n-i}u^R_{m+i}-(-1)^n u^R_{m+n-i}v^R_i\}
=(v_nu)^R_m \quad \mbox{ for }n\geq 0.\eqno{(4.5)}$$
For $m\in \Z$, we also have:
$$\begin{array}{l}
\displaystyle{(v_nu)^R_m(p_{-1}t)-p_{-1}(v_nu)_mw}
=\displaystyle{\sum_{i=0}^{m}\binom{m}{i}((v_nu)_ip)^I_{m-1-i}w }\cr
=\displaystyle{\sum_{i\in \N}\sum_{j\in \N}\binom{m}{i}\binom{n}{j}(-1)^j
(\{v_{n-j}u_{i+j}p-(-1)^nu_{n+i-j}v_jp\})^I_{m-1-1}w }\cr
=\displaystyle{\sum_{i\in \N}\sum_{j\in \N}
\sum_{h\in \N}\binom{m}{i}\binom{n}{j}(-1)^{j+h}\binom{n-j}{h}
v_{n-j-h}(u_{i+j}p)^I_{m-i-1+h}w}\cr
\mbox{}\qquad-\displaystyle{\sum_{i\in \N}\sum_{j\in \N}\sum_{h\in\N}(-1)^{n+h}
\binom{m}{i}\binom{n}{j}\binom{n-j}{h}
(u_{i+j}p)^I_{m-1+n-i-j-h}v_hw }\cr
\mbox{}\qquad+\displaystyle{\sum_{i\in \N}\sum_{j\in \N}\sum_{h\in \N}(-1)^{j+n+1+h}
\binom{m}{i}\binom{n}{j}\binom{n+i-j}{h}
u_{n+i-j-h}(v_jp)^I_{m-i-1+h}w}\cr
\mbox{}\qquad+\displaystyle{\sum_{i,j,h\in\N}(-1)^{i+h}
\binom{m}{i}\binom{n}{j}\binom{n+i-j}{h}
(v_jp)^I_{n+m-j-1-h}u_hw }\cr
=\displaystyle{\sum_{i,j,h\in \N}\binom{m}{i}\binom{n}{j}(-1)^{j+h}\binom{n-j}{h}
\{A^{n-j-h}B^{i+j}\!-\!(-1)^{n-j+1}C^{i+j}D^h\} }\cr
\mbox{}\qquad+\displaystyle{\sum_{i,j,h\in\N}\binom{m}{i}(-1)^{j+n+1+h}\binom{n}{j}\binom{n+i-j}{h}
\{E^{n+i-j-h}F^j\!-\!(-1)^{n+i-j}G^jH^h\}}\cr
=\displaystyle{(1+B)^m(A-1-B)^n-(1+C)^m(-1+D-C)^n}\cr
\mbox{}\qquad-\displaystyle{(1+E\!-\!1)^m(-E\!+\!1\!+\!F)^n+(1-1+H)^m(1-H+G)^n},  
\end{array}$$
where the monomials $A^sB^k$, $B^sD^k$, $E^sF^k$ and 
$G^sH^k$ in the last two polynomials denote \\
$\displaystyle{v_{s}(u_{k}p)^I_{m+n-s-k-1}w}, 
\displaystyle{(u_{s}p)^I_{m-1+n-s-k}v_kw}, 
\displaystyle{u_s(v_kp)^I_{m+n-k-s-1}w}, 
\mbox{ and } \displaystyle{(v_sp)^I_{n+m-s-1-k}u_kw}$, \\
respectively. On the other hand, for $n$-th normal product $\ast_n$, we have:
$$\begin{array}{l}
\displaystyle{(v \ast_n u)^R_m(p_{-1}t)-p_{-1}(v_n u)_mw }
=\!\displaystyle{\sum_{h\in\N}\binom{n}{h}(-1)^h\{v_{n-h}u_{m+h}
-(-1)^nu_{m+n-h}v_h\}(p_{-1}w)}\cr
=\!\displaystyle{\sum_{h\in\N}\binom{n}{h}(-1)^h
\left\{ v_{n-h}\binom{h\!+\!m}{j}(u_jp)^I_{m+h-1-j}w
+\sum_{j\in\N}\binom{n\!-\!h}{j}
(v_jp)^I_{n-h-1-j}u_{m+h}w \right\}}\cr 
\mbox{}\quad+\displaystyle{
\sum_{h,j\in\N}\binom{n}{h}(-1)^{n+h+1}\left\{u_{m+n-h}\binom{h}{j}(v_jp)^I_{h-1-j}w 
+\binom{m+n-h}{j}(u_jp)^I_{n+m-h-j}v_hw \right\}}\cr
=\!\displaystyle{\sum_{h,j\in\N}\binom{n}{h}(-1)^h\binom{m+h}{j}A^{n-h}B^j
+\sum_{h,j\in\N}\binom{n}{h}(-1)^{h+m+1}\binom{n+m-h}{j}G^jH^h }\cr
\mbox{}\qquad+\displaystyle{\sum_{h,j\in\N}\binom{n}{h}(-1)^{n+h+1}\binom{h}{j}E^{n+m-h}F^j
+\sum_{h,j\in\N}\binom{n}{h}(-1)^{h}\binom{n-h}{j}C^jD^{m+h}}\cr
=\!\displaystyle{(1+B)^m(A-1-B)^n+H^m(1+G-H)^n-E^m(-E+1+F)^n}\cr
\mbox{}\qquad\displaystyle{-(1+C)^m(-1-C+D)^n}.
\end{array}$$ 
Therefore, we have another Associativity 
$$(v\ast_nu)^R_m r=(v_nu)^R_m r \quad \mbox{ for }m\geq 0 \eqno{(4.6)}$$
for $r\in R$. The remaining Associativity we have to prove is the case for $m<0$ and $n<0$, 
but we will not prove it directly. Fortunately, 
the results above are enough to prove that $R$ is a $V$-module as we will see. 
Since the weight of $R$ is bounded below, 
$$\overline{V}:=\{v^R(z)=\sum_{h\in \Z} v^R_h z^{-h-1}\mid v\in V\}$$
satisfy the lower truncation property and so they are quantum operators (or 
called weak vertex operators).  
Furthermore they are mutually commutative because of (4.4) and so 
they generate a local system $\widetilde{V}$ in $\End(R)[[z,z^{-1}]]$ 
by using normal products $\ast_n$. 
Since $\omega^R(z)$ is also a Virasoro element 
of $\widetilde{V}$, 
$\widetilde{V}$ is a vertex algebra with the same grading on $\bar{V}$, that is, 
$\omega^R(z)\ast_1v^R(z)=\wt(v)v^R(z)$ by (4.5).
By viewing $n$-th normal products in $\widetilde{V}$ as 
$n$-th products in $V$, we have a vertex algebra epimorphism 
$\psi:\widetilde{V} \to V$.  Since 
$\widetilde{V}$ stabilizes 
$$0\subseteq Q\boxtimes W\subseteq R,$$ 
we have $(\Ker(\psi))^2=0$ and so $\Ker(\psi)$ is a $V$-module. 
Suppose $\widetilde{V}\not=\overline{V}$, then 
$$\CF=\{ v^R(z)\ast_n u^R(z)-(v_nu)^R(z)\mid v,u\in V \}\not=0$$  
Since $v^R(z)\ast_nu^R(z)=(v_nu)^R(z)$ for $n\geq 0$ by (4.5),  
the weights of elements in $\CF$ have a lower bound.  
Choose $0\not=\alpha(z)=v^R(z)_n\ast u^R(z)-(v_nu)^R(z)\in \CF$ so that 
$\alpha(z)$ has minimal weight among them. 
Since $\omega^R(z)\ast_k v^R(z)=(\omega_kv)^R(z)$ for $k\geq 0$, 
$\omega^R(z)\ast_k\alpha(z)$ is a linear sum of elements in $\CF$ and so 
$\omega^R(z)\ast_k(\alpha(z))=0$ for $k\geq 2$ because of the 
minimality of $\wt(\alpha(z))$. For $m\geq 0$, 
since $\alpha(z)_m=0$ on $R$ by (4.6), we have: 
$$0=[\omega^R(z)_k,\alpha(z)_m]
=(\omega_0\alpha(z))_{k+m}+k(\wt(\alpha))\alpha(z)_{k+m-1}=(-k-m+k(\wt(\alpha)))\alpha(z)_{k+m-1}$$
for any $k\in \Z$ and so $\alpha(z)_h=0$ for any $h$, 
which contradict to $\alpha(z)\not=0$. 
Therefore, $\widetilde{V}=V$ and we know that $R$ is a $V$-module. 

The remaining thing we have to prove is to define an $(z\frac{d}{dz}-L(-1)z)$-nilpotent 
intertwining operator $\CY^0$ of type $\binom{R}{P\,\,W}$. 
It is natural to expect that $\CY^0$ satisfies $\frac{d}{dz}\CY^0(q,z)=I(\omega_0q,z)$.
However, a coefficient $(L(-1)p)^I_{0}$ of $I(L(-1)q,z)$ at $z^{-1}$ may not be zero 
whereas $\frac{d}{dz}\CY^0(w,z)$ does not have $z^{-1}$-term.  
Fortunately, we have that for any $v\in V, n\in \Z$,
$$\begin{array}{l}
\displaystyle{v^R_n(L(-1)p)^I_{0}-(L(-1)p)^I_{0}v^R_n}
=\!\displaystyle{\sum_{i=0}^{\infty}\binom{n}{i}(v^P_iL(-1)p)^I_{n-i} }\cr
=\!\displaystyle{\sum_{i=0}^{\infty}\binom{n}{i}\{(L(-1)v^P_ip)^I_{n-i}+\!(iv^P_{i-1}p)^I_{n-i}\} }
=\!\displaystyle{\sum_{i=0}^{\infty}\binom{n}{i}\{(-n+i)(v^P_ip)^I_{n-i-1}+\!(iv^P_{i-1}p)^I_{n-i}\}}\cr
=\!\displaystyle{-n\sum_{i=0}^{\infty}\binom{n-1}{i}(v^P_ip)^I_{n-i-1}
+\!\sum_{i=0}^{\infty}n\binom{n-1}{i-1}(v^P_{i-1}p)^I_{n-i}}=0
\end{array}$$
and so $(L(-1)p)^I_{0}$ is a $V$-homomorphism. 
We now set 
$$\begin{array}{rll}
 p(z)w:=&\displaystyle{\left\{\int I(L(-1)p,z)-(L(-1)p)^I_{0}z^{-1})dz\right\}w
+p_{-1}w }\quad 
&\mbox{ for } w\in W\cr
=&\displaystyle{\sum_{m\in \Z,m\not=0} \frac{1}{-m}(L(-1)p)^I_{m}z^{-m}w+p_{-1}w} &\cr
p(z)r:=&\{I(p,z)-(L(-1)p)^I_{0}z^{-1}\}r  \qquad &\mbox{ for } r\in Q\boxtimes U.
\end{array} \eqno{(4.7)}$$ 
Let us show that $p(z)$ satisfies Commutativity with respect to all actions 
$a(z)$ of $a\in V$. 
By direct calculation, we have: 
$$\begin{array}{l}
[a(z), p(x)]-\dsp{\sum_{n\in \Z}
\sum_{i=0}^{\infty}\binom{n}{i}(a_ip)^I_{n-1-i}z^{-n-1}}
=\dsp{\sum_{m\not=0,n\in \Z}[a_nz^{-n-1}, (L(-1)p)^I_m
\frac{-1}{m}x^{-m}] }\cr
\mbox{}\qquad=\dsp{\sum_{m\not=0,n\in \Z}\sum_{i=0}^{\infty}\binom{n}{i}(a_iL(-1)p)^I_{m+n-i}
\frac{-1}{m}x^{-m}z^{-n-1}} \cr
\mbox{}\qquad=\dsp{\sum_{m\not=0,n\in \Z}\sum_{i=0}^{\infty}\binom{n}{i}\{(L(-1)a_ip)^I_{m+n-i}
+i(a_{i-1}p)^I_{m+n-i}\}
\frac{-1}{m}x^{-m}z^{-n-1}}\cr
\mbox{}\qquad=\dsp{\sum_{m\not=0,n\in \Z}\sum_{i=0}^{\infty}\binom{n}{i}\{(-m-n+i)(a_ip)^I_{m+n-i-1}
+i(a_{i-1}p)^I_{m+n-i}\}\frac{-1}{m}x^{-m}z^{-n-1}} \cr
\mbox{}\qquad=\dsp{\sum_{m\not=0,n\in \Z}\sum_{i=0}^{\infty}\binom{n}{i}(-m)(a_ip)^I_{m+n-i-1}
\frac{-1}{m}x^{-m}z^{-n-1}}\cr
\mbox{}\qquad=\dsp{\sum_{m\not=0,n\in \Z}
\sum_{i=0}^{\infty}\binom{n}{i}(a_ip)^I_{m+n-i-1}x^{-m}z^{-n-1}},
\end{array}$$
and so 
$$[a(z),p(x)]=\sum_{m,n\in \Z}\sum_{i\in \N}
\binom{n}{i}(a_ip)^{I}_{(m+n-i-1)}x^{-m}z^{-n-1}.$$
Since the weights of $R$ are bounded below, there is $N$ such that 
$a_ip=0$ for $i\geq N$ and so we have Commutativity:
$$
(x-z)^{N+1}[a(z), p(x)]
=(x-z)^{N+1}\dsp{(\sum_{i=0}^N\sum_{r\in \Z}(a_ip)^I_{r-i-1}\sum_{n\in \Z}\binom{n}{i}x^{n-r}z^{-n-1})}=0.$$
We then extend it by 
$$\CY^0(v_np,z)=\Res_{x}\{(x-z)^nv(x)p(z)-(-z+x)^np(z)v(x)\}$$
for $v\in V, n\in \Z$, then 
$\CY^0$ is a $(z\frac{d}{dz}-L(-1)z)$-nilpotent intertwining operator and so 
$$\CY(u,z)=\sum_{i=0}^K\frac{1}{i!}(zL(-1)-z\frac{d}{dz})^i\CY^0(u,z)\log^iz 
\quad \mbox{ for }u\in P$$
is an intertwining operator of $P$ from $W$ to $R$. 
Since $CY(u,z)=\CY_{{\rm Ker}(\rho),U}^{\boxtimes}(u,z)$ for $u\in {\rm Ker}(\rho)$, 
the space spanned by images of $\CY$ contains $Q\boxtimes U$. 
Furthermore, since $\CY(q,z)w\in p_{-1}w+Q\boxtimes U$, we get that $\CY$ is surjective. 

This completes the proof of Theorem. 
\prend

As a corollary of Lemma \ref{KOP}, we have:

\begin{thm}\label{SR}
If a simple VOA $V$ contains a rational subVOA $W$ containing 
the same Virasoro element of $V$, then $V$ is projective.   
\end{thm}

\pr 
Suppose false and let $\rho:P\to V$ be a proper covering of $V$ and choose $p\in P_0$ such that 
$\rho(p)=\1$. Viewing $P$ as a $W$-module, $P$ is a direct of irreducible $W$-modules. 
We may choose $p$ in a simple $W$-submodule $R$. Then $\rho(R)=W$ and so 
$R$ is isomorphic to $W$. Therefore $L(-1)p=0$. 
Since $P$ is a $V$-module, for $v\in V$, there is $N_v\in \N$ such that 
$v_mp=0$ for $m\geq N_v$. Then since $0=L(-1)v_mp=-mv_{m-1}p$, we have $v_kp=0$ 
for $k\in \N$, which implies $\Ker(\rho)=0$ by (4.2) and we have a contradiction. 
\prend \vspace{-4mm}

\subsection{Projective covers for semi-rigid modules}
In this subsection, we will prove the following:

\begin{thm}\label{PDS}
If $C$ is a semi-rigid irreducible module 
and $D\xrightarrow{\alpha} C$ is a covering of $C$, 
then there is a covering $R\xrightarrow{\rho} V$ such that 
$D$ is a homomorphic image of $R\boxtimes D$. In particular, 
if $V$ has a projective cover $P_V$ and $P_V\boxtimes C$ is 
finitely generated, then $C$ has a projective cover $P_C$ which is 
isomorphic to a direct summand of $P_V\boxtimes C$. 
\end{thm}

\pr
Since $C$ is semi-rigid, we have the following epimorphism:
$$C\boxtimes Q \xrightarrow{{\rm id}_C\boxtimes \epsilon} 
C\boxtimes (\widetilde{C}\boxtimes C) \xrightarrow{\mu_C} 
(C\boxtimes \widetilde{C})\boxtimes C \xrightarrow{e_C\boxtimes{\rm id}_C} C,$$
where $\epsilon:Q\to \widetilde{C}\boxtimes C$ is an embedding, 
$e_{\widetilde{C}}:\widetilde{C}\boxtimes C\to V$, 
$e_{\widetilde{C}}\epsilon:Q\to V$ is a covering and 
$\mu_C$ is a natural isomorphism for the associativity of products of intertwining 
operators, see (5.3). 
By the pull back of $D\xrightarrow{\alpha} C$, we may choose 
an isomorphism $\mu_D$ such that 
we have the following commutative diagram: 
$$\begin{array}{ccccccc}
     &  & C\boxtimes (\widetilde{C}\boxtimes D) &\xrightarrow{\mu_D}& 
(C\boxtimes \widetilde{C})\boxtimes D &\to &D \cr
      &   & \qquad\qquad \qquad \downarrow {\rm id}_C\boxtimes({\rm id}_{\widetilde{C}}\boxtimes \alpha) &  
& \qquad \qquad \downarrow {\rm id}_{C\boxtimes\widetilde{C}}\boxtimes \alpha   &   & 
\quad\downarrow \alpha\cr
C\boxtimes Q &\to& C\boxtimes (\widetilde{C}\boxtimes C) &\xrightarrow{\mu_C}& 
(C\boxtimes \widetilde{C})\boxtimes C &\to &C.
\end{array}$$
Let $R$ be a minimal submodule of $\widetilde{C}\boxtimes D$ with respect to 
inclusion such that 
$({\rm id}_{\widetilde{C}}\boxtimes \alpha)(R)=Q$. Then since 
$D\xrightarrow{\alpha} C$ is a covering and 
$(e_{W}\boxtimes {\rm id}_C)\mu_C({\rm id}_{C}\boxtimes \alpha)({\rm id}_{C}\boxtimes \epsilon_R)(C\boxtimes R)=C$, 
we have $(e_{W}\boxtimes {\rm id}_D)\mu_D({\rm id}_{C}\boxtimes \epsilon_R)(C\boxtimes R)=D$ and so 
$D$ is a homomorphic image of $C\boxtimes R$. 
By the choice of $R$, $e_{\widetilde{C}}({\rm id}_{\widetilde{C}}\boxtimes \alpha):R\to V$ 
is a covering of $V$, which proves the first statement in the theorem. 

If $V$ has a projective cover $P_V$ and $P_V\boxtimes C$ has a finite length of 
composition series, 
then any covering $D$ of $C$ is a homomorphic image of $P_V\boxtimes C$ and so 
the length of composition series of $D$ has an upper bound. Therefore, 
a covering of $C$ has a maximal one, which is projective 
by Proposition \ref{Covering} and isomorphic to a direct summand of $P_V\boxtimes C$. 

This completes the proof of Theorem \ref{PDS}. 
\prend

We easily have the following corollary. 

\begin{cry} 
If $V$ is projective as a $V$-module, then every semi-rigid irreducible module is 
projective. 
\end{cry}

\begin{thm}\label{FLT}
Let $V$ be a $C_2$-cofinite VOA of CFT-type and $C$ a semi-rigid $V$-module. Then 
$$0\rightarrow {\rm \Ker}(\rho_C)\boxtimes W 
\xrightarrow{\epsilon\boxtimes {\rm id}_W} 
P_C\boxtimes W 
\xrightarrow{\rho_C\boxtimes {\rm id}_W} C
\boxtimes W\rightarrow 0  \eqno{(4.8)}$$
is exact, where $0\rightarrow {\rm Ker}(\rho_C)\xrightarrow{\epsilon} P_C
\xrightarrow{\rho_C} C\to 0$ is a projective cover of $C$. 
\end{thm}

\pr
Since $V$ is $C_2$-cofinite, all fusion products are finitely generated and 
the products of intertwining operators satisfy Associativity. 
Let $\rho_V:P_V\to V$ be a projective cover of $V$. 
By Theorem 14, there is a surjection $\alpha:P_V\boxtimes C\rightarrow P_C$. 
Set $J=\Ker(\alpha)$. Since $C\in {\rm mod}(V)$, 
we have the following commutative exact diagram by (4.1):
$$\begin{array}{ccccccccc}
0 &\rightarrow&     J      &\rightarrow&      J      &\rightarrow&    0       &&  \cr
      \dto&  &\dto&  &\dto& &\dto& & \cr
0 &\rightarrow &\Ker(\rho_V)\boxtimes C &\rightarrow& P_V\boxtimes C  &
\xrightarrow{\rho_V\boxtimes {\rm id}_C}& V\boxtimes C=C &\rightarrow&  0 \cr
      \dto&  &\dto&  &\quad \dto \alpha& &\dto& &\dto\cr
0 &\rightarrow&  \Ker(\rho_C)      &\xrightarrow{\epsilon}&      D      
&\xrightarrow{\rho_C}&    C       &\rightarrow&  0 \cr
     &  &\dto&  &\dto& &\dto& &\cr
 &&     0      &\rightarrow&      0      &\rightarrow&    0       &&   
\end{array}$$ 
From Theorem \ref{FPP}, we also have an exact sequence 
$$  0\rightarrow \Ker(\rho_V)\boxtimes (C\boxtimes W) \rightarrow  
P_V\boxtimes (C\boxtimes W) \rightarrow C\boxtimes W\rightarrow 0.$$
Since the fusion products satisfy the associativity and preserve the 
right exactness, we have an exact sequence:
$$ 0\rightarrow (\Ker(\rho_V)\boxtimes C)\boxtimes W \rightarrow 
(P_V\boxtimes C)\boxtimes W \rightarrow C\boxtimes W\rightarrow 0. $$
Again, since $\boxtimes W$ preserves the right exactness, 
we have the commutative exact diagram:
$$\begin{array}{ccccccccc}
0 &\rightarrow&    J\boxtimes W    &\rightarrow&      J\boxtimes W    
&\rightarrow&    0       &&  \cr
    \dto  &  &\sigma_1 \dto\quad&  &\sigma_2\dto\quad& &\dto& & \cr
0 &\rightarrow &(\Ker(\rho_V)\boxtimes C)\boxtimes W &\rightarrow& 
(P_V\boxtimes C)\boxtimes W  &\rightarrow
  & C\boxtimes W &\rightarrow&  0 \cr
  &  &\dto&  &\dto& &\dto& &\dto\cr
  & &     \Ker(\rho_C)\boxtimes W      &\xrightarrow{\epsilon\boxtimes {\rm id}_W}&      
D\boxtimes W &\rightarrow
  &    C\boxtimes W       &\rightarrow&  0 \cr
  &  &\dto&  &\dto& &\dto& &\cr
  & &     0    &\rightarrow&      0      &\rightarrow&    0,       &&  
\end{array}$$
which implies that $\epsilon\boxtimes {\rm id}_W$ is injective. 

This completes the proof of Theorem \ref{FLT}. 
\prend

Since $V\boxtimes W\cong W$, the following easily comes from Theorem \ref{PDS}. 

\begin{cry}\label{AP}
Let $V$ be a simple $C_2$-cofinite VOA of CFT-type and 
assume that all $V$-modules are semi-rigid. 
If $V$ is projective as a $V$-module, then all $V$-modules are 
projective. In particular, $V$ is rational.
\end{cry}

\section{Products of intertwining operators} 
\subsection{Fusion products}
In this section, we assume that the desired products of 
intertwining operators satisfy the associativity and 
the desired fusion products are well-defined as $V$-modules. 
The aim in this section is to use homomorphisms and fusion product intertwining operators 
to present products of intertwining operators. 

We first recall the analytic part on the 
composition of intertwining operators (with logarithmic terms) from \cite{HD}. 
From now on, let $\{A,B,C,D,E,F\}$ be a set of $C_1$-cofinite $V$-modules.  We choose 
$a\in A, b\in B, c\in C, d'\in D'$. As Huang showed,  
for intertwining operators $\CY_1\in \CI\binom{D}{A\,\,E}$, 
$\CY_2\in \CI\binom{E}{B\,\,C}$, 
$\CY_3\in \CI\binom{D}{F\,\,C}$ and $\CY_4\in \CI\binom{F}{A\,\,B}$, 
the formal power series (with logarithmic terms) 
$$ \begin{array}{ccc}
\langle  d',\CY_1(a,x)\CY_2(b,y)c\rangle  & \mbox{ and } &
\langle  d',\CY_3(\CY_4(a,x-y)b,y)c\rangle 
\end{array}$$
are absolutely convergent in  $\Delta_1=\{(x,y)\in \C^2\mid|x|>|y|>0\}$ 
and $\Delta_2=\{(x,y)\in \C^2 \mid |y|>|x-y|>0\}$, respectively, 
and can all be analytically extended to multi-valued analytic functions on 
$$M^2=\{(x,y)\in \C^2 \mid xy(x-y)\not=0 \}.$$
As he did, we are able to lift them to single-valued analytic functions 
$$\begin{array}{ccc}
E(\langle  d,\CY_1(a,x)\CY_2(b,y)c\rangle ) &\mbox{ and }&
E(\langle  d,\CY_3(\CY_4(a,x-y)b,y)c\rangle ) 
\end{array}\eqno{(5.1)}$$
on the universal covering $\widetilde{M^2}$ of $M^2$. 
Single-valued liftings are not unique as he remarked, but 
the existence of such functions is enough for our arguments. 
The important fact is that if we fix $A,B,C,D$, then these functions 
are given as solutions of the same differential equations. 
Therefore, for $\CY_1\in \CI\binom{D}{A\,\,E},
\CY_2\in \CI\binom{E}{B\,\,C}$ there are $\CY_5\in \CI\binom{D}{A\boxtimes B\,\,C}$ 
and 
$\CY_6\in \CI\binom{D}{B\,\,A\boxtimes C}$ 
such that 
$$\begin{array}{cl}
E(\langle d',\CY_1(a,x)\CY_2(b,y)c\rangle)
=E(\langle d',\CY_5(\CY_{A,B}^{\boxtimes}(a,x-y)b,y)c\rangle) &\mbox{ and}\cr
E(\langle d',\CY_2(\CY_4(a,x-y)b,y)c\rangle)
=E(\langle d',\CY_6(a,x)\CY_{B,C}^{\boxtimes}(b,y)c\rangle). &
\end{array}\eqno{(5.2)}$$ 
We note that the right hand sides of (5.2) are usually expressed as linear sums, say, 
$$E(\langle d',\CY_1(a,x)\CY_2(b,y)c\rangle)=\sum_i E(\langle d',\CY_{1i}(\CY_{2i}(a,x-y)b,y)c\rangle).  $$
However, for each $i$, from the maximality of fusion products, 
there is a homomorphism $\xi_i\in \Hom_V(A\boxtimes B, \Image(\CY_{2i}))$ such that 
$\CY_{2i}=\xi_i\CY_{A,B}^{\boxtimes}$. Since $\CY_5:=\sum_i \CY_{1i}\xi_i$ is 
an intertwining operator in 
$\CI\binom{D}{A\boxtimes B\,\,C}$, we get the expressions (5.2).  For example,  
a natural isomorphism $\mu: (A\boxtimes B)\boxtimes C \to 
A\boxtimes(B\boxtimes C)$ is 
given by 
$$E(\langle d',\mu
\CY_{A\boxtimes B, C}^{\boxtimes}(\CY_{A,B}^{\boxtimes}(a,x-y)b,y)c\rangle)
=E(\langle d',\CY_{A,B\boxtimes C}^{\boxtimes}(a,x)\CY_{B,C}^{\boxtimes}(b,y)c\rangle). 
\eqno{(5.3)}$$
In order to simplify the notation, 
we often use the notations $\CY_1\CY_2$ and $\CY_3(\CY_4)$ to distinguish two types 
$\CY_1(\ast,x)\CY_2(\ast, y)$ and $\CY_3(\CY_4(\ast,x)\ast,y)$, respectively. 

\subsection{Skew symmetric and adjoint intertwining operators}
In his paper \cite{HV}, Huang has explicitly defined a skew symmetry intertwining operator 
$\sigma_{12}(\CY)\in \CI\binom{C}{B,A}$ and an adjoint intertwining operator 
$\sigma_{23}(\CY)\in \CI\binom{B'}{A\,\,C'}$ for $\CY\in \CI\binom{C}{A\,\,B}$ when 
$\CY$ has no logarithmic terms. Even if $\CY\in \CI\binom{C}{B\,\,A}$ has logarithmic terms, 
there are isomorphisms $\sigma_{12}:\CI\binom{C}{A\,\,B} \to \CI\binom{C}{B\,\,A}$ 
and 
$\sigma_{23}:\CI\binom{C}{A\,\,B} \to \CI\binom{B'}{A\,\,C'}$ as follows. 
By considering a path $\{z=\frac{1}{2}e^{\pi i t}x \mid 0\leq t\leq 1\}$, 
there is $\widetilde{\CY}\in \CI\binom{C}{A\,\,B}$ such that 
$$E(\langle c', \widetilde{\CY}(a,z)\sigma_{12}(Y^B)(b,x)\1\rangle)
=E(\langle c',\CY(b,x)\sigma_{12}(Y^A)(a,z)\1\rangle). \eqno{(5.4)}$$
Rewriting them, we have that 
$$\begin{array}{rl}
\mbox{the left side of (5.4)}=&E(\langle c', \widetilde{\CY}(a,z)e^{L(-1)x}b \rangle)
=E(\langle c', e^{L(-1)x}\widetilde{\CY}(a,z-x)b \rangle) \cr
=&E(\langle e^{L(1)x}c', \widetilde{\CY}(a,z-x)b \rangle)  \qquad \mbox{ and}\cr
\mbox{the right side of (5.4)}
=&E(\langle c', \CY(b,x)e^{L(-1)(z)}a \rangle)
=E(\langle c', e^{L(-1)z}\CY(b,x-z)a \rangle) \cr
=&E(\langle e^{L(1)x}c', e^{L(-1)(z-x)}\CY(b,x-z)a \rangle).
\end{array}$$
Since $\langle e^{L(1)x}c', \widetilde{\CY}(a,z-x)b \rangle$ 
and $\langle e^{L(1)x}c', e^{L(-1)z}\CY(b,x-z)a \rangle$ are multivalued 
rational functions on $\{(x,z)\mid x\not=z\}$, 
we may choose $\sigma_{12}$ so that $\sigma_{12}(\CY)=\widetilde{\CY}$, that is, 
$$\sigma_{12}(\CY)(a,z-x)b=e^{L(-1)(z-x)}\CY(b,x-z)a. \eqno{(5.5)}$$  
Similarly, for $\CY\in \CI\binom{C'}{A\,\,B}$ and 
canonical intertwining operators $\CY_{C,C'}^{V'}$ and $\CY_{B',B}^{V'}$ 
induced from inner products, there is $\CY_4\in \CI\binom{B'}{A\,\,C}$ such that 
$$ E(\langle \1, \CY_{C,C'}^{V'}(c,x)\CY(a,y)b\rangle)=
E(\langle \1, \CY_{B',B}^{V'}(e^{L(-1)(x-y)}\CY_4(a,y-x)c,y)b\rangle).\eqno{(5.6)}$$
Therefore, a map $\CY\to \CY_4$ gives an isomorphism $\sigma_{23}:
\CI\binom{C}{A\,\,B}\cong \CI\binom{B'}{A\,\,C'}$. 

We don't get an explicit formula because what we will need is just 
an existence of isomorphism $\sigma_{23}(\CY)$. 
Set $\sigma_{123}=\sigma_{12}\sigma_{23}$. 

In (5.2), we used $\CY^{\boxtimes}$ as the second intertwining operator  
of products. 
Not only the second one, we can also use it for the first 
one at the same time. Actually, for 
$\CY_5(\CY_{A,B}^{\boxtimes})$ with $\CY_5\in \CI\binom{D}{A\boxtimes B\,\,C}$,  
we have $\sigma_{23}^{-1}\sigma_{12}^{-1}(\CY_5)\in 
\CI\binom{(A\boxtimes B)'}{C\,\,D'}$ and so 
there is $\delta\in \Hom_V(C\boxtimes D',(A\boxtimes B)')$ 
such that $\sigma_{23}^{-1}\sigma_{12}^{-1}(\CY_5)=\delta \CY_{C,D'}^{\boxtimes}$. 
Therefore we have: 

\begin{prn} 
For any $\CY_1\in \CI\binom{D}{A\,\,E}, 
\CY_2\in \CI\binom{E}{B\,\,C}$ and $\CY_3\in \CI\binom{E}{A\,\,B}, 
\CY_4\in \CI\binom{D}{E\,\,C}$, 
there are $\xi\in \Hom_V(B\boxtimes C, (A\boxtimes D')')$ and 
$\delta\in \Hom_V(B\boxtimes C, (C\boxtimes D')')$ such that 
$$\begin{array}{l}
\langle d',\CY_1(a,z_1)\CY_2(b,z_2)c\rangle
=\langle d',\sigma_{23}(\CY_{A,D'}^{\boxtimes})(a,z_1)
\xi \CY_{B,C}^{\boxtimes}(b,z_2)c\rangle, \quad\mbox{ and} \cr
\langle d',\CY_4(\CY_3(a,z_1)b,z_2)c\rangle
=\langle d',\sigma_{123}(\CY_{C,D'}^{\boxtimes})
(\delta \CY_{A,B}^{\boxtimes}(a,z_1)b,z_2)c\rangle, 
\end{array} \eqno{(5.7)}$$
for any $a\in A, b\in B, c\in C, d'\in D'$.  
In other words, the space spanned by the products $\CY_1\CY_2$ of 
$\CY_1\in \CI\binom{D}{A\,\,\ast}$ and $\CY_2\in \CI\binom{\ast}{B\,\,C}$ is 
isomorphism to $\Hom_V(B\boxtimes C,(A\boxtimes D')')$ and 
the space spanned by the product $\CY_3(\CY_4)$ of 
$\CY_3\in \CI\binom{\ast}{A\,\,B}$ and $\CY_4\in \CI\binom{D}{\ast\,\,C}$ is 
isomorphism to $\Hom_V(A\boxtimes B,(C\boxtimes D')')$. 

\end{prn} \vspace{-4mm}

\subsection{Semi-rigidity and intertwining operators}
The aim in this subsection is to describe the semi-rigidity in terms of intertwining operators. 
We assume that $V$ is simple. 
For a $V$-module $U$, let ${\rm rad}^V(U)$ denote the intersection of maximal submodules $M$ 
with $U/M\cong V$. Since $V$ is simple, $U/{\rm rad}^V(U)$ is a direct sum of copies of $V$. 
Let $W$ be a simple $C_1$-cofinite $V$-module and assume $W$ is not semi-rigid. 
By the definition of semi-rigid, for any $\widetilde{W}\in {\rm mod}(V)$ satisfying (1) in Definition 1, 
$e_W\in \Hom_V(W\boxtimes \widetilde{W},V)$ and a maximal submodule $Q$ of $\widetilde{W}\boxtimes W$ 
containing ${\rm rad}^V(\widetilde{W}\boxtimes W)$, 
$(e_W\boxtimes {\rm id}_W)\mu({\rm id}_W\boxtimes \epsilon)$ is not surjective, in other words, 
$$(e_W\boxtimes {\rm id}_W)\mu({\rm id}_W\boxtimes \epsilon)(W\boxtimes Q) \subseteq 
 \Ker(e_W\boxtimes {\rm id}_W).$$ 
Since $\mu$ is isomorphism and $W$ is simple, the above implies that 
$$
\mu(({\rm id}_W\boxtimes \epsilon)(W\boxtimes {\rm rad}^V(\widetilde{W}\boxtimes W)))
+
\Ker(e_W\boxtimes {\rm id}_W)=
(W\boxtimes \widetilde{W})\boxtimes W, \eqno{(5.8)}$$ 
for any $e_{W}:W\boxtimes \widetilde{W}\to V$, where 
$\epsilon:{\rm rad}^V(\widetilde{W}\boxtimes W)\subseteq \widetilde{W}\boxtimes W$
is an embedding. 

We then rewrite (1.1) by products of intertwining operator as follows:
Since $e_{\widetilde{W}}\CY_{\widetilde{W},W}^{\boxtimes}\in \CI\binom{V}{\widetilde{W},W}$, 
from (5.2) and (5.7), there are $\CY_5$ and $\delta\in \Hom_V(W\boxtimes \widetilde{W}, (W\boxtimes W')')$ 
such that  
$$\begin{array}{l}
E(\langle a', \sigma_{12}(Y^{W})(w,x)
e_{\widetilde{W}}\CY_{\widetilde{W},W}^{\boxtimes}(\widetilde{w},y)w^1
\rangle)=E(\langle a',\CY_5(\CY_{W,\widetilde{W}}^{\boxtimes}(w,x-y)\widetilde{w},y)w^1 \cr
\mbox{}\hspace{5cm}=E(\langle a', \sigma_{123}(\CY_{W,W'}^{\boxtimes})
(\delta \CY_{W,\widetilde{W}}^{\boxtimes}(w,x-y)\widetilde{w},y)
w^1 \rangle) 
\end{array}\eqno{(5.9)}$$
for any $w, w^1\in W$, $\widetilde{w}\in \widetilde{W}$, and $a'\in W'$. Therefore, we have:

\begin{lmm}\label{NSR} Let $V$ be a simple VOA. 
If a simple $V$-module $W$ is not semi-rigid, then 
the image of $\delta$ does not have a factor isomorphic to $V$ and 
$\Ker(\delta) +{\rm rad}^V(W\boxtimes \widetilde{W})=W\boxtimes \widetilde{W}$ 
for any $e_{\widetilde{W}}$. 
\end{lmm}

\section{Orbifold model}
In this section, we will consider an orbifold model. 
Let $T$ be a VOA and $g$ an automorphism of $T$ order $n$ and set 
$T=\oplus_{h=0}^{n-1} T^{(h)}$ with $T^{(h)}:=\{t\in T\mid g(t)=e^{2\pi \sqrt{-1}h/n}t\}$. 
Then $V:=T^{(0)}=T^g$ is a subVOA. Let $D$ be a $V$-module and we assume that    

\noindent  
(A1)\quad  every $T^{(i)}\boxtimes_V D$ is a $V$-module for $i=0,...,n-1$ and \\
(A2)\quad  every $T^{(i)}$ is $C_1$-cofinite as a $T^{(0)}$-module. \\

The necessary condition we need is an associativity of fusion products, that is, 
$T^{(i)}\boxtimes_V (T^{(j)}\boxtimes_V D)\cong (T^{(i)}\boxtimes_V T^{(j)})\boxtimes_V D$. 

Set $W^{(i)}=T^{(i)}\boxtimes_V D$ and $W=W^{(0)}\oplus \cdots \oplus W^{(n-1)}$. 
We note $W^{(0)}=V_V\boxtimes D=D$.  
Since all $T^{(i)}$ are $C_1$-cofinite as $V$-modules, 
there is $\CY\in \CI\binom{W}{T\,\,W}$ such that 
$$ E(\langle w', \CY(t,x)\CY_{T,D}^{\boxtimes}(t^1,y)d\rangle) 
=E(\langle w', \CY_{T,D}^{\boxtimes}(Y(t,x-y)t^1,y)d\rangle) \eqno{(6.1)}$$
for $t,t^1\in T$, $w'\in W'$ and $d\in D$ by (5.2). 
We note $W^{(i)}=T^{(i)}\boxtimes D=\coprod_r (W^{(i)}_r)$ with 
formal weight spaces $W^{(i)}_r:=(T^{(i)}\boxtimes D)_r$ and so we can define its 
restricted dual $(W^{(i)})'=\coprod_r \Hom(W^{(i)}_r,\C)$.    
From the definition of $\CY$ and Commutativity of vertex operators $Y$ of $T$, 
we have that for any $t^1,t^2,t^3\in T$, 
$$ \begin{array}{rl}
E(\langle w',\CY(t^1,x)\CY(t^2,y)\CY_{T,D}^{\boxtimes}(t^3,z)d\rangle) 
=&\!E(\langle w',\CY(t^1,x)\CY_{T,D}^{\boxtimes}(Y(t^2,y-z)t^3,z)d\rangle)\cr 
=&\!E(\langle w',\CY_{T,D}^{\boxtimes}(Y(t^1,x-z)Y(t^2,y-z)t^3,z)d\rangle)\cr
=&\!E(\langle w',\CY_{T,D}^{\boxtimes}(Y(t^2,y-z)Y(t^1,x-z)t^3,z)d\rangle)\cr
=&\!E(\langle w',\CY(t^2,y)\CY(t^1,x)\CY_{T,D}^{\boxtimes}(t^3,z)d\rangle),
\end{array}$$
which implies Commutativity of $\{\CY(t,z) \mid t\in T\}$ on $W$. We also have 
$$ \begin{array}{rl}
E(\langle w',\CY(t^1,x)\CY(t^2,y)\CY_{T,D}^{\boxtimes}(t^3,z)p\rangle) 
=&\!E(\langle w',\CY_{T,D}^{\boxtimes}(Y(t^1,x-z)Y(t^2,y-z)t^3,z)p\rangle)\cr
=&\!E(\langle w',\CY_{T,D}^{\boxtimes}(Y(Y(t^1,x-y)t^2,y-z)t^3,z)p\rangle)\cr
=&\!E(\langle w',\CY(Y(t^1,x-y)t^2,z_2)\CY_{T,D}^{\boxtimes}(t^3,z)p\rangle).
\end{array}$$
Furthermore, taking $t^1=\1$ in (6.1), we have 
$$\begin{array}{rl}
E(\langle w',\CY(t,x)d\rangle)=&E(\langle w',\CY(t,x)\CY_{T,D}^{\boxtimes}(\1,y)d\rangle)
=E(\langle w',\CY_{T,D}^{\boxtimes}(Y(t,x-y)\1,y)d\rangle)\cr
=&E(\langle w',\CY_{T,D}^{\boxtimes}(e^{(x-y)L(-1)}t,y)p\rangle) 
=E(\langle w',\CY_{T,D}^{\boxtimes}(t,y+x-y)d\rangle).
\end{array}$$
By setting $y=\frac{2}{3}x$, we obtain 
$\CY(t,x)d=\CY_{T,D}^{\boxtimes}(t,x)d$ for $t\in T, d\in D$. 

Furthermore, we have:

\begin{prn}\label{OSC}
Let $T$ be a VOA and $g\in {\rm Aut}(T)$ of order $n$. 
If all $T^{(k)}$ are $C_1$-cofinite as $T^{(0)}$-modules and $T^{(k)}\boxtimes T^{(1)}$ are  
$T^{(0)}$-modules,   
then all $T^{(k)}$ are simple currents.
\end{prn}

\pr 
We will prove that $T^{(1)}$ is a simple current. 
The proofs for the others are similar. 
Take $T^{(1)}$ as $D$ in the above argument and we 
assert that $U^{(0)}:=T^{(n-1)}\boxtimes T^{(1)}$ is irreducible. 
Suppose false and let $B$ be a proper submodule of $U^{(0)}$ and set 
$$E=\langle \mbox{ the coefficients of }\CY(t^{(1)},y)b 
\mid t^{(1)}\in T^{(1)}, b\in B\rangle_{\C}. $$ 
Since $E$ is a submodule of $T^{(1)}$ and 
$$\CY(t^{(n-1)},x)\CY(t^{(1)},y)b=\CY(Y(t^{(n-1)},x-y)t^{(1)},y)b, \eqno{(6.2)}$$
for any $t^{(1)}\in T^{(1)}$, $t^{(n-1)}\in T^{(n-1)}$ and $b\in B$,  
there are $t\in T^{(1)}$ and $b\in B$ such that $\CY(t,y)b\not=0$ and so 
we have $E\not=0$. Therefore, we have $E=T^{(1)}$ since $T^{(1)}$ is simple. 
Since the coefficients of $\CY(t^{(1)},y)b$ spans $T^{(1)}$, 
the left hand side of (6.2) implies that the coefficient of (6.2) spans 
$T^{(n-1)}\boxtimes T^{(1)}$. 
On the other hand, the right hand side of (6.2) implies that 
those spans $B$ since 
the coefficients of $Y(t^{(n-1)},x-y)t^{(1)}$ are actions of $V$. 
Therefore, $T^{(n-1)}\boxtimes T^{(1)}=V$. 
For any $V$-module $W$, $W=V\boxtimes W=(T^{(n-1)}\boxtimes T^{(1)})\boxtimes W
\cong T^{(n-1)}\boxtimes(T^{(1)}\boxtimes W)$, which implies that 
$T^{(1)}\boxtimes W$ is simple.   
\prend

As another application, we have:

\begin{prn} 
Under the same assumption as in Proposition \ref{OSC}, 
if $0\to B\to D\xrightarrow{\rho}V$ is a non-split extension of $V$ by a simple $V$-module $B$, then 
$T\boxtimes_VD$ is a $T$-module.  In particular, if $T$ is projective as a $T$-module, then $V$ is projective 
as a $V$-module. 
\end{prn}

\pr 
By the previous arguments, it is sufficient to show that 
$\CY(t,x)$ is a formal integer power series. 
Since $(T^{(n-1)}\boxtimes T^{(1)})\boxtimes D\cong D$ and 
$T^{(1)}\boxtimes D$ is an indecomposable module with a homomorphic image 
$T^{(1)}\boxtimes V$, all weights of $T^{(1)}\boxtimes D$ are integers. 
Hence the formal powers are all integer powers. 
So the remaining thing is to show that $\CY(t,x)$ does not have $\log z$-terms. 
If $L(0)$ acts on $D$ semisimply, then so does on $T^{(k)}\boxtimes D$ and 
every $\CY_{T^{(k)},D}^{\boxtimes}(t,z)$ is a formal integer power series, which 
implies that $T\boxtimes_VD$ is a $T$-module. 
If $L(0)$ does not act on $D$ semisimply, then $L(0)^{nil}D=B$ and so $B\cong V$ as $V$-modules. 
Therefore, $L(0)^{nil}:(T\boxtimes_VD)/(T\boxtimes_VB)\to T\boxtimes_VB\cong T$ is an 
isomorphism as $T$-modules. 
If $\CY_{T,D}^{\boxtimes}(t,z)p$ is not a formal integer power series for $t\in T$ and $p\in D$, 
then we have: 
$$\CY_{T,D}^{\boxtimes}(t,z)p=\sum_{m\in \Z} t_{(0,m)}pz^{-m-1}
+\sum_{m\in \Z} t_{(1,m)}pz^{-m-1}\log z.$$
Set $\CY^{(0)}(t,z)=\sum_{m\in \Z} t_{(0,m)}z^{-m-1}$ and 
$\CY^{(1)}(t,z)=\sum_{m\in\Z} t_{(1,m)}z^{-m-1}\in 
\Hom((T\boxtimes_V D)/(T\boxtimes_VB),T\boxtimes_V B)[[z,z^{-1}]]$. 
Since $\dim \CI\binom{T}{T\,\,T}=1$ as $T$-modules, we may assume 
$\CY^{(0)}(t,z)d \equiv Y(t,z)d \mod{T\boxtimes_VB[[z,z^{-1}]]}$ and 
$\CY^{(1)}(t,z)a=Y(t,z)\overline{a}$, 
where $\overline{a}$ denotes ${\rm id}_T\boxtimes \rho(a)$ and $a\in T\boxtimes_VD$.  
Therefore, we have $L(0)^{nil}\CY^{(0)}(t,z)a=Y(t,z)\bar{a}$ and 
$Y(t,z)\CY^{(0)}(t,z)L(0)^{nil}a=Y(t,z)\bar{a}$. In other words, $L(0)^{nil}$ commutes with 
$\CY^{(0)}$. 
By (2.1), we have 
$$\begin{array}{rl}
w_{(1,n)}d=&-(L(0)^{nil})t_{(0,n)}d+(L(0)^{nil}t)_{(0,n)}d+t_{(0,n)}(L(0)^{nil}d) \cr
  =&-(L(0)^{nil})t_{(0,n)}d+t_{(0,n)}(L(0)^{nil}d)=0. 
\end{array}$$
Therefore, $\CY(t,z)$ is a formal integer power series and so 
$T\boxtimes_V D$ is a $T$-module. 
\prend

\section{Geometrically modified module}
In this section, we will explain the theory of composition-invertible power series 
and their actions on modules for the Virasoro algebra developed in \cite{HD} and 
then we will extend them for logarithmic intertwining operators.  
From now on, $q_x$ denotes $e^{2\pi ix}$ for variables $x\not=\tau$ 
to distinguish it from $q=e^{2\pi i\tau}$.  
Let $A_j$  $(j=1,2,...)$ be the complex numbers defined by 
$$\frac{1}{2\pi i}(q_y-1)=\left( \exp\left(-\sum_{j=1}^{\infty} A_j y^{j+1}\frac{\pd}{\pd y}\right) \right)y $$
and set 
$$\CU(q_x)={q_x}^{L(0)}(2\pi i)^{L(0)}e^{-\sum_{j=1}^{\infty} A_jL(j)}.$$ 
Clearly, $\CU(q_x)={q_x}^{L(0)}\CU(1)$. The important operator 
is $\CU(1)$, which satisfies 
$$ \CU(1)\CY(w,x)\CU(1)^{-1}=\CY(\CU(q_x)w, q_x-1)=\CY({q_x}^{L(0)}\CU(1)w,q_x-1)=\CY[\CU(1)w,x]  \eqno{(7.1)}$$
for an intertwining operator $\CY$, see \cite{Zh} for $\CY[\cdot,x]$.

\subsection{Trace functions}
We first consider $q$-traces of geometrically-modified module 
operators with one more variable $z$:
$$\Psi_U(v;z,\tau):=\Tr_U Y(\CU(q_z)v, q_z)q^{(L(0)-c/24)} \eqno{(7.2)}$$
for a $V$-module $U$ and $v\in V$, 
where $\Tr_U$ denotes a trace on $U$ and $c$ is the central charge of $V$. We note that for an ordinary 
trace function, we can consider the trace functions for 
not only actions of $V$ but also actions of a $V$-module $W$ on $U$ by $\CY\in \CI\binom{U}{W\,\,U}$. 
Namely, we can define a trace function 
$$\Psi_U(\CY;w;z,\tau):
=\Tr_U (\CY(\CU(q_z)w, q_z)q^{(L(0)-c/24)}) \qquad w\in W. \eqno{(7.3)}$$

We have to note that $L(0)$ may not be semisimple on a $V$-module $U$. 
Then we will understand $q^{L(0)}$ on $U$ as 
$$q^{L(0)}:=q^{(\wt+L(0)^{nil})}
=q^{\wt}(e^{2\pi i \tau L(0)^{nil}})=
q^{\wt}\sum_{j=0}^{\infty}\frac{(2\pi i\tau L(0)^{nil})^j}{j!}.$$   
We note that since $L(0)^{nil}$ is a nilpotent operator and commutes with all grade-preserving 
operators, there are no terms of form $q^{r}\tau^j$ in ordinary trace functions for $j>0$.  

We note that when we consider a trace function of a simple module $W$ on a $V$-module $U$, 
as we explained at the end of \S 2.2, we can ignore the $\log z$ terms of 
$\CY_{W,U}^U\in \CI\binom{U}{W\,\,U}$ and so the necessary grade-preserving operator of $w\in W_r$ in 
$\CY_{W,U}^U(w,z)=\sum_i\sum_m w_{(i,m)}^{\CY}z^{-m-1}\log^iz$ is $w^{\CY}_{(0,r-1)}$. 
Therefore, by setting $\CU(1)w=\sum_{r} w^r$ with homogeneous elements $w^r\in W_r$, we have 
$$\begin{array}{rl}
\Tr_U \CY_{W,U}^U(\CU(q_z)w, q_z)q^{(L(0)-c/24)}=&\sum_r \Tr_U 
{q_z}^{(\wt(w^{r}))}w^{r}_{(0,r-1)}q^{(-r)}q^{(L(0)-c/24)}\cr
=&\sum_r \Tr_U w^{r}_{(0,r-1)} q^{(L(0)-c/24)}.  
\end{array}\eqno{(7.4)}$$ 
Thus, (7.4) is independent of $z$. Moreover, 
it has shown in \cite{HD} that these $q$-traces are absolutely 
convergent when $0<|q|<1$ and can be analytically extended to analytic functions of 
$\tau$ in the upper-half plane.

We next consider $q$-traces of products of two geometrically-modified intertwining operators: 
$$\begin{array}{l}
\Tr_U \CY_1(\CU(q_y)\CY_{W,\widetilde{W}}^{\boxtimes}(w,x-y)\widetilde{w}, q_y)q^{(L(0)-c/24)}\cr
\Tr_U \CY_2(\CU(q_x)w,q_x)
\CY_{\widetilde{W},U}^{\boxtimes}(\CU(q_y)\widetilde{w},q_y)q^{(L(0)-c/24)}
\end{array}\eqno{(7.5)}$$ 
for $w\in W, \widetilde{w}\in \widetilde{W}$, 
$\CY_1\in \CI\binom{U}{W\boxtimes \widetilde{W}\,\,U}$, and 
$\CY_2\in \CI\binom{U}{W\,\,\widetilde{W}\boxtimes U}$. 
As we explained, the first function in (7.5) depends on $x-y$, but not on $y$. 
These formal power series (with log-terms) are 
absolutely convergent in $\Omega_1=\{(x,y,\tau)\in \C^2\oplus \CH \mid 0<|q_x-q_y|<|q_y|\}$ 
and $\Omega_2=\{(x,y,\tau)\in \C^2\oplus \CH\mid 0<|q|<|q_y|<|q_x|<1\}$, 
respectively, as shown in \cite{HD}, 
where $\CH=\{\tau\in \C\mid {\rm Im}(\tau)>0\}$ is the upper half plane. 
We extend these function analytically to multivalued 
analytic functions on  
$$M_1^2=\{(x,y,\tau)\in \C^2\times \CH \mid x\not=y+p\tau+q \quad\mbox{ for all }
p,q\in \Z\}.$$
We can lift them to 
single valued analytic functions 
$$\begin{array}{l}
\Psi_U(\CY_1(\CY_{W,\widetilde{W}}^{\boxtimes}):w,\widetilde{w};x,y,\tau):=E(\Tr_U 
\CY_1(\CU(q_y)\CY_{W,\widetilde{W}}^{\boxtimes}(w,x-y)\widetilde{w}, q_y)q^{(L(0)-c/24)}) \cr
\Psi_U(\CY_2\CY_{\widetilde{W},U}^{\boxtimes}:w,\widetilde{w};x,y,\tau):=E(\Tr_U \CY_2(\CU(q_x)w,q_x)
\CY_{\widetilde{W},U}^{\boxtimes}(\CU(q_y)\widetilde{w},q_y)q^{(L(0)-c/24)})
\end{array}\eqno{(7.6)}$$ 
on the universal covering $\widetilde{M^2_1}$. 
We need to extend a result in \cite{HV} 
for logarithmic intertwining operators.

\begin{lmm}\label{GEQ}
For every intertwining operator 
$\CY\in \CI\binom{T}{B,U}$, $w\in W$ and $b\in B$, we have 
$$\begin{array}{l}  e^{\tau L(0)}\CY(b,z)u=\CY(e^{\tau L(0)}b,e^{\tau}z)e^{\tau L(0)}u  \qquad 
\cr
q^{L(0)}\CY(\CU(q_y)b,q_y)=\CY(\CU(q_yq)b,q_yq)q^{L(0)}\qquad \mbox{ and} \cr
\CY_1(\CY_2(\CU(q_y)b,q_y-q_x)\CU(q_x)w,q_x)
=\CY_1(\CU(q_x)\CY_2(b,y-x)w,q_x).
\end{array}$$
\end{lmm}

\pr
Set $\CY(b,z)=\sum_{h=0}^K\sum_{n\in \C} b_{(h,n)}z^{-n-1}\log^h z$ and 
$y=\log z$. From \\
$\CY(L(-1)b,z)=\frac{d}{dz}\CY(b,z)$, 
we have $(L(-1)b)_{(h,n+1)}=(-n-1)b_{(h,n)}+(h+1)b_{(h+1,n-1)}$ and 
$$[L(0),b_{(h,n)}]u=(L(-1)b)_{(h,n+1)}+(L(0)b)_{(h,n)}\!=\!
(-n-1)b_{(h,n)}\!+\!(h+1)b_{(h+1,n)}\!+\!(L(0)b)_{(h,n)}.$$
Using the notations $(\alpha\otimes \beta)b_{(k,n)}u=(\alpha b)_{(k,n)}\beta u$, we have:
$$\begin{array}{l}
L(0)(b_{(h,n)}u)=(-n\!-\!1\!+L(0)\otimes 1+1\otimes L(0))b_{(h,n)}u+(h+1)b_{(h+1,n)}
\qquad \mbox{ and}\cr
L(0)^m(b_{(h,n)}u)=\sum_{j=0}^m\binom{m}{j}
(-n\!-\!1\!+L(0)\otimes 1+1\otimes L(0))^{m-j}(h+1)\cdots(h+j)b_{(h+j,n)}u.
\end{array}$$
for $m\geq 1$, where $(h+1)\cdots (h+j)=1$ for $j=0$. Using these notation, we obtain:
$$\begin{array}{l}
e^{\tau L(0)}\CY(b,z)u
=\displaystyle \sum_ne^{\tau L(0)}(\sum_{h=0}^K b_{(h,n)}uy^h)e^{(-n-1)y}\cr
=\!\!\displaystyle{\sum_n\sum_{m=0}^{\infty}
\sum_{h=0}^{\infty} \frac{L(0)^m\tau^m}{m!}(\sum_h b_{(h,n)}uy^he^{(-n-1)y})}\cr
=\!\!\displaystyle{\sum_n\sum_{m,h}\frac{\tau^m}{m!}\sum_{j=0}^{m}\!\binom{m}{j}
(L(0)\!\otimes 1\!+\!1\otimes L(0)\!-\!n\!-\!1)^{m-j}(h\!+\!1)\cdots(h\!+\!j)
b_{(h+j,n)}uy^he^{(-n-1)y}} \cr
\mbox{By replacing $h+j$ and $m-j$ by $k$ and $i$, 
respectively, $e^{\tau L(0)}\CY(b,z)u$ equals}\cr
\displaystyle{\sum_n\sum_k\sum_{i=0}^{\infty}\sum_{j=0}^{\infty}
\frac{\tau^{i}(L(0)\otimes 1+1\otimes L(0)\!-\!n\!-\!1)^{i}}{i!}\frac{1}{j!}
(k-j+1)\cdots(k)\tau^jy^{k-j}b_{(k,n)}ue^{(-n-1)y}}\cr
=\!\!\displaystyle{\sum_n\sum_{k=0}^{\infty}e^{\tau(L(0)\otimes 1
+1\otimes L(0)-(n+1))}b_{(k,n)}u(y+\tau)^ke^{(-n-1)y}}\cr 
=\!\!\displaystyle{\sum_n\sum_k 
(e^{\tau L(0)}b)_{(k,n)}(e^{\tau L(0)}u) e^{(-n-1)(y+\tau)}(y+\tau)^k
=\CY(e^{\tau L(0)}b, e^{\tau+y})e^{\tau L(0)}u}\cr
=\!\!\CY(e^{\tau L(0)}b, e^{\tau}z)e^{\tau L(0)}u,  
\end{array}$$
which proves the first equation. Replacing $\tau$ and $y$ by $2\pi i\tau$ and $2\pi i y$, respectively, 
we have the second equation.
The third comes from $\CU(1)\CY(b,x)=\CY(\CU(q_x)b,q_x-1)\CU(1)$ and the second 
equation. 
\prend \vspace{-4mm}

\section{Transformations}
Let $V$ be a simple $C_2$-cofinite VOA of CFT-type. 
We fix two irreducible $V$-modules $W$ and $\widetilde{W}$ such that 
$\Hom_V(\widetilde{W}\boxtimes W,V)\not=0$.  
In this section, we always assume that the desired fusion products are $V$-modules and 
the products of intertwining operators of $W$ and $\widetilde{W}$ satisfy the associativity. 

\subsection{Three transformations}
We define actions $S$, $\alpha_t$, $\beta_t$ on $R_2^1$ 
for $0\leq t\leq 1$ by 
$$\begin{array}{rl}
\alpha_t&: (x,y,\tau)\mapsto (x,y+t\tau,\tau) \cr
\beta_t&: (x,y,\tau)\mapsto (x,y+t,\tau) \cr
S&: (x,y,\tau)\mapsto (-x/\tau,-y/\tau,-1/\tau) 
\end{array}. $$
In particular, we have the following commutative diagram:
$$ \begin{array}{ccc}
(x,y,\tau) &\xrightarrow{S} & (-x/\tau,-y/\tau, -1/\tau) \cr
\downarrow \beta_t & & \downarrow \alpha_t \cr
(x,y+t,\tau) &\xrightarrow{S} & (-x/\tau, -y/\tau-t/\tau,-1/\tau). 
\end{array}\eqno{(8.1)}$$

A trace function of products 
$\CY_1(\CY)$ with $\CY_1\in \CI\binom{U}{E\,\,U}, \CY\in 
\CI\binom{E}{\widetilde{W}\,\,W}$ on $U$ is 
$$\Psi_U(\CY_1(\CY):\widetilde{w},w;x,y,\tau)
=E\left(\Tr_U \CY_1\left(\CU(q_y)\CY(\widetilde{w},x-y)w,q_y\right)q^{(L(0)-c/24)}\right),
\eqno{(8.2)} $$
for $w\in W,\widetilde{w}\in \widetilde{W}$. A modular transformation $S:\tau\to -1/\tau$ on $\Psi_U$ is 
defined by 
$$
S\left(\Psi_U\right)\left(\CY_1(\CY):\widetilde{w},w;x,y,\tau\right)=\Psi_U\left(\CY_1(\CY): 
\left( \frac{-1}{\tau}\right)^{L(0)}\widetilde{w}, \left(\frac{-1}{\tau}\right)^{L(0)}w;
\frac{-x}{\tau},\frac{-y}{\tau}; \frac{-1}{\tau}\right).$$
In order to simplify the arguments, we will deal only $W$ and $\widetilde{W}$ with integer weights 
when we consider the transformation $S$. 
When $\CY_1$ is a vertex operator $Y^U$ on a module, the space spanned by trace functions 
has some modular invariance property as the author has shown in \cite{M1}. 
In particular, since we assume Condition II, that is, there are $\lambda_U\in \C$ for $U\in {\rm Irr}(V)$ such that  
$$S\left(\Psi_V\right)\left(Y^U(\CY):\widetilde{w},w;x,y,\tau\right)
=\sum_{U\in {\rm Irr}(V)} \lambda_U \left(\Psi_U\right)\left(Y^U(\CY):\widetilde{w},w;x,y,\tau\right), 
\eqno{(8.4)}$$
where ${\rm Irr}(V)$ denotes the set of irreducible $V$-modules. 
We note that $\lambda_U$ does not depend on the choice $\CY\in \CI\binom{V}{W\,\,\tilde{W}}$, 
but on $V$. 

Along a line $\CL=\{(x,y+t,\tau)\mid t\in [0,1]\}$ from 
$(x,y,\tau)$ to $(x,y+1,\tau)$, we define 
$$
\alpha_t(\Psi_U)(\CY:\widetilde{w},w;x,y,\tau):
=\Psi_U(\CY:\widetilde{w},w;x,y+t,\tau)
\eqno{(8.5)}$$
and set $\alpha=\alpha_1$. 
Since $(x,y,\tau)\to (x,y+t,\tau)$ preserves 
$\Omega_2=\{(x,y,\tau)\in \C^2\oplus H\mid |q|<|q_y|<|q_x|<1\}$, we have 
$$\begin{array}{rl}
\alpha_t(\Psi_U)(\CY_1(\CY_2):\widetilde{w},w;x,y,\tau)
=&\alpha_t(\Tr_U \CY_3(\CU(q_x)\widetilde{w},q_x)
\CY_{\widetilde{W},U}^{\boxtimes}(\CU(q_y)w,q_y)q^{(L(0)-c/24)})\cr
=&\Tr_U \CY_3(\CU(q_x)\widetilde{w},q_x)
\CY_{\widetilde{W},U}^{\boxtimes}(\CU(q_yq^t)w,q_yq^t)q^{(L(0)-c/24)}\cr
=&\Tr_U \CY_4(\CU(q_y)\CY_5(\widetilde{w},x-y)w, q_y)q^{(L(0)-c/24)}\cr
=&\Psi_U(\CY_4(\CY_5):\widetilde{w},w;x,y,\tau)
\end{array}\eqno{(8.6)}$$
for some $\CY_3$ and $\CY_4(\CY_5)$, because $\CY_{W,U}^{\boxtimes}(\CU(q_yq^t)w,q_yq^{t})$ 
is a linear combination of geometrically modified intertwining 
operators in $\CI\binom{\boxtimes}{W\,\,U}$. 

An important case is where $U=V$ and $\CY_1(\CY_2)=Y(\CY)$ 
with $\CY\in \CI\binom{V}{\widetilde{W}\,\,W}$. 
Then since $W\boxtimes V=W$ is irreducible, 
$$\alpha(\Psi_V)(Y(\CY))=e^{2\pi i\wt(W)}\Psi_V(Y(\CY)).$$ 
We set $\kappa=e^{2\pi i\wt(W)}$.  
We then define $\beta_t$ and $\beta=\beta_1$ according to a line $S^{-1}(\CL)$ by 
$$\beta_t(\Psi_U)(\CY_1(\CY_2):\widetilde{w},w;x,y,\tau)
=\Psi_U(\CY_1(\CY_2):\widetilde{w},w;x,y+t\tau,\tau) \mbox{ for any }\Psi_U.  \eqno{(8.7)}$$

\begin{prn}  
$$\beta_t (S(\Psi_V)))=S(\alpha_t(\Psi_V)).  \eqno{(8.8)}$$ 
\end{prn}

By (8.4), we will consider the following diagram: 
$$\begin{array}{cccc}
    \Psi_V(Y(\CY)) & \xrightarrow{\alpha} &  &\kappa \Psi_V(Y(\CY)) \cr
      \downarrow \quad S        &   & & \downarrow \quad S  \cr
    \sum \lambda_U\Psi_U(Y^U(\CY)) &\xrightarrow{\beta}  & 
\sum \lambda_U\beta(\Psi_U(Y^U(\CY))&\kappa\sum \lambda_U\Psi_U(Y^U(\CY)) \cr
    \end{array}$$
Therefore $\sum \lambda_U\beta(\Psi_U(Y^U(\CY))=\kappa\sum \lambda_U\Psi_U(Y^U(\CY))$.

\subsection{The image of $\beta$}

We first calculate the image of $\Psi_U(\CY^1(\CY^2):\widetilde{w},w;x,y,\tau)$ by $\beta$ 
for any $\CY^2\in \CI\binom{S}{\widetilde{W}\,\,W}$ 
and $\CY^1\in \CI\binom{U}{S\,\, U}$. 
Set $A=(W\boxtimes U)$, then we have:
$$\begin{array}{l}
\beta(\Psi_U)(\CY^1(\CY^2):\widetilde{w},w;x,y,\tau) \cr 
\mbox{}=E(\Tr_U \CY^1(\CU(e^{2\pi i(y+\tau)})
\CY^2(\widetilde{w},x-(y+\tau))w, 
e^{2\pi i(y+\tau)})q^{(L(0)-\frac{c}{24})})
\quad \cr
\mbox{}=E(\Tr_U \CY^1(\CY^2(\CU(q_x)\widetilde{w},q_x-e^{2\pi i(y+\tau)})
\CU(e^{2\pi i(y+\tau)})w, 
q_yq)q^{(L(0)-\frac{c}{24})})\quad \mbox{by Lemma \ref{GEQ}}\cr
\mbox{}=E(\Tr_{U} \sigma_{23}(\CY_{\widetilde{W},U'}^{\boxtimes})(\CU(q_x)\widetilde{w},q_x)\xi_U
\CY_{W,U}^{\boxtimes}
(\CU(q_yq)w,q_yq)q^{(L(0)-\frac{c}{24})})\cr
\mbox{}\qquad \qquad \mbox{for some }\xi_U\in \Hom_V(W\boxtimes U,(\widetilde{W}\boxtimes U')') \cr
\mbox{}=E(\Tr_U \sigma_{23}(\CY_{\widetilde{W},U'}^{\boxtimes})(\CU(q_x)\widetilde{w},q_x)
q^{(L(0)-\frac{c}{24})}
\xi_U\CY_{W,U}^{\boxtimes}(\CU(q_y)w,q_y))
\qquad \mbox{by Lemma \ref{GEQ}}\cr
\mbox{}=E(\Tr_U \sigma_{23}(\CY_{\widetilde{W},U'}^{\boxtimes})(\CU(q_x)\widetilde{w},q_x)
\xi_U q^{(L(0)-\frac{c}{24})}
\CY_{W,U}^{\boxtimes}(\CU(q_y)w,q_y))\cr
\mbox{}=E(\Tr_A \CY_{W,U}^{\boxtimes}(\CU(q_y)w,q_y)
\sigma_{23}(\CY_{\widetilde{W},U'}^{\boxtimes})(\CU(q_x)\widetilde{w},q_x)\xi_U q^{(L(0)-c/24)})\cr
\mbox{}\qquad \qquad \mbox{because the trace is symmetric}\cr
\mbox{}=E(\Tr_A \sigma_{123}(\CY_{A,A'}^{\boxtimes}(\delta_U \CY_{W,\widetilde{W}}^{\boxtimes})(\CU(q_y)w,
q_y-q_x)\CU(q_x)\widetilde{w},q_x)q^{(L(0)-c/24)}) \cr
\mbox{}\qquad \qquad \mbox{for some $\delta_U\in \Hom_V(W\boxtimes \widetilde{W},(A\boxtimes A')')$} \cr
\mbox{}=E(\Tr_A \sigma_{123}(\CY_{A,A'}^{\boxtimes})(\delta_U 
\CU(q_x)\CY_{W,\widetilde{W}}^{\boxtimes}(w,y-x)\widetilde{w},q_x)q^{(L(0)-c/24)}) 
\qquad\mbox{ by Lemma \ref{GEQ}}\cr
\mbox{}=E(\Tr_A \sigma_{123}(\CY_{A,A'}^{\boxtimes})(\delta_U 
{q_x}^{L(0)}\CU(1)e^{L(-1)(y-x)}\sigma_{12}(\CY_{W,\widetilde{W}}^{\boxtimes})
(\widetilde{w}, x-y)w,q_x)q^{(L(0)-c/24)}) \cr
\mbox{}\qquad \qquad \mbox{by skew symmetry intertwining operator, see (5.5)}. 
\end{array}$$
Set $L[-1]=L(-1)+L(0)$ (see \cite{Zh}).   
Then we get $U(1)e^{L(-1)z}=e^{(2\pi i)L[-1]z}U(1)$ from (7.1) and 
the above equals the following:
$$\begin{array}{l}
\mbox{}=E(\Tr_A \sigma_{123}(\CY_{A,A'}^{\boxtimes})(\delta_U 
{q_x}^{L(0)}e^{(2\pi i)L[-1](y-x)}\CU(1)\sigma_{12}(\CY_{W,\widetilde{W}}^{\boxtimes})
(\widetilde{w},x-y)w,q_x)q^{(L(0)-c/24)}).
\end{array}$$
As we explained in the paragraph after (7.5), the pair of terms ${q_x}^{L(0)}$ and $q_x$ in the 
above expression has no influence and the next term is $e^{(2\pi i)L[-1](y-x)}$. However, since 
$o_0(L[-1]u)=0$ for any $u\in \widetilde{W}\boxtimes W$, we finally have 
$$\begin{array}{l}
\beta(\Psi_U)(Y^U(\CY_{\widetilde{W},W}^V):\widetilde{w},w;x,y,\tau) \cr 
\mbox{}\quad=E(\Tr_A \sigma_{123}(\CY_{A,A'}^{\boxtimes})
\CU(q_x)\delta_U\sigma_{12}(\CY_{W,\widetilde{W}}^{\boxtimes}))
(\widetilde{w}, x-y)w,q_x)q^{(L(0)-c/24)}). 
\end{array}\eqno{(8.9)}$$
In particular, we have the following lemma. \\

\noindent
\begin{lmm}\label{SOT}
$\beta(\Psi_{U})(\CY^1(\CY^2))$ is again 
an ordinary trace function of an intertwining operator of $W\boxtimes \widetilde{W}$ on 
some module.
\end{lmm}

Therefore, we have:

\begin{lmm}\label{DLT}
When $\CY^1=Y^U$ and $\CY^2=e_{\widetilde{W}}\CY_{\widetilde{W}\,\,W}^{\boxtimes}$, 
$\delta_U$ coincides with $\delta$ in (5.9). In particular, $\delta_U$ does not depend on $U$. 
\end{lmm}

\pr
We express the definitions of $\xi_U$ and $\delta_U$ in a short way:
$$
Y^U(\CY_{\widetilde{W},W}^V)
=\sigma_{23}(\CY_{\widetilde{W},U'}^{\boxtimes})\xi_U\CY_{W,U}^{\boxtimes}
\quad\mbox{ and }\quad 
\CY_{W,U}^{\boxtimes}\sigma_{23}(\CY_{\widetilde{W},U'}^{\boxtimes})\xi_U=
\sigma_{123}(\CY_{A,A'}^{\boxtimes})(\delta_U \CY_{W,\widetilde{W}}^{\boxtimes}).
\eqno{(8.10)}$$ 

For $a'\in A'$, $\widetilde{w}\in \widetilde{W}$, $w,w^1\in 
W$ and $u\in U$, let us consider 
$$\begin{array}{l}
\langle a', \CY_{W,U}^{\boxtimes}(w^1,x)
Y^U(e_{\widetilde{W}}\CY_{\widetilde{W},W}^{\boxtimes}(\widetilde{w},y-z)w,z)u\rangle 
\end{array} \eqno{(8.11)}$$
into two ways. 
Set $B={\rm Image}(\delta_U)$, then 
there is $\CY_{B,W}^{(U\boxtimes A')'}$ such that  
$$\begin{array}{rl}
\mbox{(8.11)}&=
\langle a', \CY_{W,U}^{\boxtimes}(w^1,x)
\sigma_{23}(\CY_{\widetilde{W},U'}^{\boxtimes})(\widetilde{w},y)
\xi_U\CY_{W,U}^{\boxtimes}(w,z)u\rangle \cr
&=
\langle a', \sigma_{123}(\CY_{A,A'}^{\boxtimes})
(\delta_U\CY_{W,\widetilde{W}}^{\boxtimes}(w^1,x-y)\widetilde{w},y)
\CY_{W,U}^{\boxtimes}(w,z)u\rangle \cr
&=
\langle a', \sigma_{123}(\CY_{U,A'}^{\boxtimes})
\CY_{B,W}^{(U\boxtimes A')'}(\delta_U\CY_{W,\widetilde{W}}^{\boxtimes}(w^1,x-y)\widetilde{w},y-z)
w,z)u\rangle. 
\end{array}$$
On the other hand, there is $\CY_{W,V}^{W}
\in \CI\binom{W}{W\,\,V}$ and $\epsilon\in \Hom_V(W,(U\boxtimes A')')$ such that 
$$\begin{array}{rl}
\mbox{(8.11)}=& 
\langle a', \CY_{W,U}^{\boxtimes}(Y_{W,V}^{W}
(w^1,x-z)
e_{\widetilde{W}}\CY_{\widetilde{W},W}^{\boxtimes}(\widetilde{w},y-z)w,z)u\rangle \cr
=&\langle a', \sigma_{123}(\CY_{U,A'}^{\boxtimes})(\epsilon Y_{W,V}^{W}
(w^1,x-z)
e_{\widetilde{W}}\CY_{\widetilde{W},W}^{\boxtimes}(\widetilde{w},y-z)w,z)u\rangle  
\end{array}$$
for any $a'\in A'$ and $u\in U$. We note $\CY_{W,V}^W\in \C\sigma_{12}(Y^W)$. 
Therefore, we have 
$$
\epsilon Y_{W,V}^{W}
(w^1,x-z)
e_{\widetilde{W}}\CY_{\widetilde{W},W}^{\boxtimes}(\widetilde{w},y-z)w
=\CY_{B,W}^{(U\boxtimes A')'}
(\delta_U\CY_{W,\widetilde{W}}^{\boxtimes}(w^1,x-z)\widetilde{w},y-z)w.$$
Since the image of $\epsilon$ is $W$, we obtain 
$$
\epsilon Y_{W,V}^{W}
(w^1,x-z)
\CY_{\widetilde{W},W}^{\boxtimes}(\widetilde{w},y-z)w
=\CY_{B,W}^W
(\delta_U\CY_{W,\widetilde{W}}^{\boxtimes}(w^1,x-z)\widetilde{w},y-z)
w$$
for some $\CY_{B,W}^W$. 
Thus, $\delta_U$ essentially coincides with $\delta$ in (5.9), which does not 
depend on the choice of $U$. 
\prend

\section{Key Theorem}
\begin{thm}\label{NCE} 
If a simple $C_2$-cofinite $V$ of CFT-type satisfies Condition I and Condition II, 
then all simple $V$-modules are semi-rigid. 
Furthermore, if $V$ is a rational VOA of CFT-type satisfying Condition I, 
then $\Psi_U$ has nonzero coefficient in $S(\Psi_V)$ for every simple $V$-module $U$.
\end{thm}

\pr 
Let $W$ be an irreducible module. As we showed in Lemma \ref{DLT}, 
$$\begin{array}{rl}
 \beta(\sum \lambda_U\Psi_U)(Y(e_{\widetilde{W}}\CY_{\widetilde{W},W}^{\boxtimes}))
 =&\sum \lambda_U\beta(\Psi_U)(Y(e_{\widetilde{W}}\CY_{\widetilde{W},W}^{\boxtimes})) \cr
 =&\sum \lambda_U\Psi_{W\boxtimes U}
(\CY_{B,U}^U(\delta \CY_{W,\widetilde{W}}^{\boxtimes})),
\end{array}$$
where $B={\rm Image}(\delta)$ and $U$ runs over ${\rm Irr}(V)$. 
On the other hand, since $\beta(S(\Psi_V))=S(\alpha(\Psi_V))$, we obtain 
$$\beta(\sum \lambda_U\Psi_U(Y(e_{\widetilde{W}}\CY_{\widetilde{W},W}^{\boxtimes}))
=\kappa(\sum \lambda_U\Psi_U(Y(e_{\widetilde{W}}\CY_{\widetilde{W},W}^{\boxtimes})),$$
where $\kappa=e^{2\pi i\wt(W)}$. Therefore, we have 
$$\sum \lambda_U\Psi_{W\boxtimes U}
(\CY_{B,U}^U(\delta \CY_{W,\widetilde{W}}^{\boxtimes}))
=\kappa(\sum \lambda_U\Psi_U(Y(e_{\widetilde{W}}\CY_{\widetilde{W},W}^{\boxtimes})).$$
Suppose that $W$ is not semi-rigid. 
Since $V$ is $C_2$-cofinite, all weights of modules $U$ are rational numbers. 
For any natural integer $n$, $V$ satisfies Condition I and II 
if and only if $V^{\otimes n}$ does and $\Psi_U$ appears in $S(\Psi_V)$ if and only if $\Psi_{U^{\otimes n}}$ 
appears in $S(\Psi_{V^{\otimes n}})$ and $W$ is semi-rigid if and only if so is $W^{\otimes n}$. 
Therefore, by taking $V^{\otimes n}$ and $W^{\otimes n}$ as $V$ and $W$ if necessary, we may assume 
that the conformal weight $\wt(W)$ of $W$ is an integer.  
By Lemma \ref{NSR}, $\Ker(\delta)+\Ker(e_{\widetilde{W}})=\widetilde{W}\boxtimes W$. 
Set $Q=\Ker(\delta)\cap \Ker(e_{\widetilde{W}})$ and $W\boxtimes \widetilde{W}/Q=Q^1\oplus Q^2$ with 
$Q^1=\Ker(e_{\widetilde{W}})/Q$ and $Q^2=\Ker(\delta)/Q\cong V$. 
Then $\Psi_{W\boxtimes U}
(\CY_{B,U}^U(\delta \CY_{W,\widetilde{W}}^{\boxtimes}))$ are all given by 
traces on $Q^1$ and $\Psi_U(Y(e_{\widetilde{W}}\CY_{\widetilde{W},W}^{\boxtimes})$ are 
all given by traces on $Q^2$. We hence have 
$$
 \sum \lambda_U\Psi_{W\boxtimes U}
(\CY_{B,U}^U(\delta \CY_{W,\widetilde{W}}^{\boxtimes}))=0, $$
which contradicts to 
$\sum \lambda_U\Psi_U(Y(e_{\widetilde{W}}\CY_{\widetilde{W},W}^{\boxtimes})
\not=0$. 
Therefore, $W$ is semi-rigid. Since $W$ is arbitrary, all irreducible $V$-modules 
are semi-rigid.  

We next prove the second statement and assume that $V$ is rational. 
We first show $\lambda_{V'}\not=0$. Choose a simple module $U$ so that 
$\lambda_U\not=0$. Set $W=U'$ and consider the trace function of 
$e_{\widetilde{W}}\CY_{\widetilde{W},W}^{\boxtimes}$ in 
$\beta(\Psi_U)(e_{\widetilde{W}}\CY_{\widetilde{W},W}^{\boxtimes})$. 
It has a nonzero scalar times of 
$$\Psi_{W\boxtimes U}(e_{\widetilde{W}}\CY_{\widetilde{W},W}^{\boxtimes})$$
and so it has a term 
$\Psi_{V'}(e_{\widetilde{W}}\CY_{\widetilde{W},W}^{\boxtimes})$ with 
a nonzero coefficient.  
On the other hand, for any $V$-modules $R\not=U$, 
$\beta(\Psi_R(e_{\widetilde{W}}\CY_{\widetilde{W},W}^{\boxtimes})$ has no entries 
of $\Psi_{V'}(e_{\widetilde{W}}\CY_{\widetilde{W},W}^{\boxtimes})$. 
Therefore, $\Psi_{V'}(e_{\widetilde{W}}\CY_{\widetilde{W},W}^{\boxtimes})$ has 
nonzero coefficient in 
$\beta(\sum \lambda_U
\Psi_U(Y(e_{\widetilde{W}}\CY_{\widetilde{W},W}^{\boxtimes}))$.

The remaining thing is to prove $\lambda_U\not=0$ for 
every simple module $U$. 
Set $W=U'$.
As we showed, $\lambda_{V'}\not=0$ and so 
there is a simple $V$-module $R$ with $\lambda_R\not=0$ such that 
$\beta(\Psi_R)(Y^R(\CY_{\widetilde{W},W}^V))$ has nonzero coefficient 
at $\Psi_{V'}(Y^U(\CY_{\widetilde{W},W}^V))$. Then since 
$\Hom_V(R\boxtimes W, V')\not=0$, we have $R=(W)'=U$ and so $\lambda_U\not=0$ as 
we desired. \\
This completes the proof of Theorem \ref{NCE}. 
\prend \vspace{-4mm}

\section{Rationality of orbifold model}
In this section, we will show that an orbifold model $R^{(0)}$ satisfies Condition II and prove the 
rationality of $T^{(0)}$ under the assumption that $T^{(0)}$ is $C_2$-cofinite and $T$ is rational. 

\begin{thm}\label{NST}
Let $g$ be a finite automorphism of a VOA $T$ of order $n$ and assume that 
a fixed point subVOA $V:=T^g$ is a simple $C_2$-cofinite of CFT-type and satisfies Condition I, 
If $T$ satisfies Condition II, then so does $V$. 
In particular, if $T$ is rational and $V$ satisfies Condition I, 
then $V$ is rational and every simple $V$-module is a submodule of some $g^k$-twisted $T$-module. 
\end{thm}

\pr	
Set $V=T^{(0)}$. 
As Dong, Li and Mason have shown in \cite{DLM}, 
there is a $g$-twisted simple $T$-module $M$, 
say $M=\oplus_{m=0}^{\infty}M_{\lambda+m/n}$. 
Then for each $i=0,\ldots,n-1$, 
$W^{(i)}=\oplus_{m=0}^{\infty}M_{\lambda+m+i/n}$ is a simple $V$-module and 
we may assume that $T^{(j)}\boxtimes W^{(k)}=W^{(k+j)}$ since 
$T^{(j)}$ is a simple current by Proposition 20. From Condition I, 
there are an irreducible $V$-module $\widetilde{W}^{(0)}$ and 
a surjection $e_{\widetilde{W}^{(0)}}:\widetilde{W}^{(0)}\boxtimes W^{(0)}\to V$. 
Set $\CY=e_{\widetilde{W}^{(0)}}\CY_{\widetilde{W}^{(0)}\,\, W^{(0)}}^{\boxtimes}$ and 
consider a trace function $\Psi_T(Y(\CY))$ of $Y(\CY)$ on $T$. 

We first consider the images of $\Psi_T(Y(\CY))$ by $\alpha^k$. 
Since $W^{(0)}\boxtimes T^{(j)}=W^{(j)}$, we have
$\alpha^k(\Psi_{T^{(j)}})(Y(\CY))=e^{2k\pi i\wt(W^{(j)})}\Psi_{T^{(j)}}(Y(\CY))$. 
Therefore, we have:
$$\alpha^k(\Psi_T(Y(\CY)))=\alpha^k(\sum_{j=0}^{n-1} \Psi_{T^{(j)}}(Y(\CY)))
=\sum_{j=0}^{n-1}e^{2k\pi i\lambda} e^{2\pi i kj/n}\Psi_{T^{(j)}}(Y(\CY)),$$ 
which coincides with $e^{2k\pi i\lambda}$-times of the trace function of $g^k$ 
$$\Tr_T o(g^k) Y^T(\CU(q_y)
e_{\widetilde{W}}\CY_{\widetilde{W}^{(0)}\,\, 
W^{(0)}}^{\boxtimes}(\widetilde{w},x-y)w, q_y)q^{(L(0)-c/24)}$$
on $T$ for each $k=0,\cdots,n-1$. 
Therefore, 
$$\Psi_V(Y(\CY))=
\frac{1}{n}\sum_{k=0}^{n-1}e^{-2k\pi i\lambda}\alpha^k\Psi_{T}(Y(\CY)).$$
On the other hand, 
since $T$ is rational and $C_2$-cofinite, 
$S(\Psi_T)$ is a linear combination of trace functions 
$\Psi_U$ on $T$-modules $U$. In particular, 
it is a linear combination of trace functions on $V$-modules.   
By Lemma \ref{SOT}, $\beta(S(\Psi_T))$ is a linear combination 
of trace functions of intertwining operators. Iterating it, 
$\beta^k(S(\Psi_T))$ are all linear combinations of trace functions of 
intertwining operators. 
Since $S\alpha^k=\beta^k S$, 
$S\alpha^k (\Tr_T Y^T(\CY^k))\CU(q_y)$ are all linear combinations of 
trace functions and so does 
$$S(\Psi_V(Y(\CY)))=\frac{1}{n}\sum_{k=0}^{n-1}e^{-2k\pi i\lambda} 
S(\alpha^k(\Psi_{T}(Y(\CY)))).$$
The second statement comes from Theorem \ref{NCE} and \ref{NST}. 
\prend 

Since every rational $C_2$-cofinite VOA of CFT-type satisfies Condition II, we have:

\begin{cry}\label{OBF}
Let $g$ be a finite automorphism of a VOA $T$ and assume that 
a fixed point subVOA $V:=T^g$ is a $C_2$-cofinite of CFT-type and satisfies Condition I, 
If $T$ is rational and $V$ satisfies Condition I, 
then $V$ is rational and every simple $V$-module is a submodule of some $g^k$-twisted $T$-module. 
\end{cry}

\end{document}